\documentclass[11pt]{article}
\usepackage{graphics,amsmath,amssymb}
\usepackage{latexsym}
\usepackage{epsfig}
\usepackage{hyperref}
\usepackage{xcolor}

\def\a{\alpha}
\def\be{\beta}

\def\epsilon{\varepsilon}
\def\ga{\gamma}

\def\la{\lambda}
\def\La{\Lambda}
\def\phi{\varphi}
\def\si{\sigma}

\def\om{\omega}

\newtheorem{theorem}{Theorem}[section]
\newtheorem{lemma}[theorem]{Lemma}

\newtheorem{corollary}[theorem]{Corollary}

\newtheorem{proposition}[theorem]{Proposition}
\newtheorem{remark}[theorem]{Remark}

\def\Z{{\mathbb Z}}
\def\N{{\mathbb N}}
\def\C{{\mathbb C}}
\def\R{{\mathbb R}}
\def\d{\,\hbox{\rm d}}

\newenvironment{Proof}{\removelastskip\par\medskip
\noindent{\em Proof.} \rm}{\penalty-20\null$\square$\par\medbreak}

\newenvironment{Proofy}{\removelastskip\par\medskip
\noindent{\em Proof of Theorem~\ref{th:inv.ingham}.} \rm}{\penalty-20\null$\square$\par\medbreak}

\newenvironment{Proofv}{\removelastskip\par\medskip
\noindent{\em Proof of Theorem~\ref{th:dir.ingham}.} \rm}{\penalty-20\null$\square$\par\medbreak}

\newenvironment{Proof1}{\removelastskip\par\medskip
\noindent{\em Proof of Theorem~\ref{th:diringham}.} \rm}{\penalty-20\null$\square$\par\medbreak}

\newenvironment{Proof2}{\removelastskip\par\medskip
\noindent{\em Proof of Proposition~\ref{pr:haraux-inv}.} \rm}{\penalty-20\null$\square$\par\medbreak}

\newenvironment{Proof3}{\removelastskip\par\medskip
\noindent{\em Proof of Theorem~\ref{th:inv.ingham1}.} \rm}{\penalty-20\null$\square$\par\medbreak}

\title{\bf Reachability problems 
for a  wave-wave  
\\
system with a memory term }

\author{Paola Loreti
\thanks{Dipartimento di Scienze di Base e Applicate per l'Ingegneria
Sezione di Matematica,
Sapienza Universit\`a di Roma,
Via Antonio Scarpa 16, 00161 Roma (Italy); e-mail: $<$loreti@dmmm.uniroma1.it$>$ }
\and Daniela Sforza
\thanks{Dipartimento di Scienze di Base e Applicate per l'Ingegneria
Sezione di Matematica, Sapienza Universit\`a di Roma,
Via Antonio Scarpa 16, 00161 Roma (Italy); e-mail: $<$sforza@dmmm.uniroma1.it$>$ }}

\begin{document}

\maketitle

\begin{abstract}
We solve the reachability problem for a coupled wave-wave system with an integro-differential term.
The control functions act on one side of the boundary. The estimates on the time is given in terms of  the parameters of the problem and they are explicitly computed thanks to Ingham type results. Nevertheless some restrictions appear in our main results. The Hilbert Uniqueness Method is briefly recalled. Our findings can be applied to concrete examples in viscoelasticity theory.
\end{abstract}


\
\section{Introduction}
The linear viscoelasticity theory has been extensively studied by many authors, that proposed several mathematical models based on experimental data to tackle such subject. A possible approach relies on the following physical assumption: the present stress is given by a functional of the past history of the deformation gradient. 
Such functionals can be represented by means of  convolution integrals.  This leads to wave equations in which a so-called memory term also appears, see the seminal papers of Dafermos \cite{D1,D2} and  \cite{RHN,LPC}. In this framework an important issue is to identify suitable class of integral kernels that match with the physical models. For example,
decreasing exponential kernels arise  in the analysis of Maxwell fluids  or Poynting -Thomson solids, see  e.g.
\cite{Pruss,Re1}. It is also noteworthy to mention  that such kernels satisfy the principle of  fading memory, {\em the memory of a simple material fades in time},
introduced in \cite{CN}. 
  
Our aim, justified by the previous remarks,
 is to investigate the reachability for  a system constituted of a  wave equation with a memory term and another wave equation coupled by lower order terms.
Precisely, given  $a\,,b\in\R$ we consider the following system
\begin{equation}\label{eq:problem-uI}
\begin{cases}
\displaystyle 
u_{1tt}(t,x) -u_{1xx}(t,x)+\beta\int_0^t\ e^{-\eta(t-s)} u_{xx}(s,x)ds+au_2(t,x)= 0\,,
\quad (0<\beta<\eta)
\\
\hskip8cm
t\in (0,T)\,,\quad x\in(0,\pi)
\\
\displaystyle
u_{2tt}(t,x) -u_{2xx}(t,x)+bu_1(t,x)= 0
\,,
\end{cases}
\end{equation}
subject to the boundary conditions
\begin{equation}\label{eq:bound-u2I}
u_1(t,0)=u_2(t,0)=0\,,\quad u_1(t,\pi)=g_1(t)\,,\quad u_2(t,\pi)=g_2(t)\qquad t\in (0,T) \,,
\end{equation}
and with null initial conditions 
\begin{equation}\label{eq:initcI}
u_i(0,x)=u_{it}(0,x)=0\qquad  x\in(0,\pi),\quad i=1,2\,.
\end{equation} 
We wish to solve a reachability problem for \eqref{eq:problem-uI} of the following type: given $T>0$ and taking
$
(u_{i}^{0},u_{i}^{1})
$, $i=1,2$,
whose regularity we will specify later, one has to
find $g_i\in L^2(0,T)$, $i=1,2$ such that the weak
solution $u$ of problem  \eqref{eq:problem-uI}-\eqref{eq:initcI}
satisfies the final conditions
\begin{equation}\label{eq:problem-u1I}
u_i(T,x)=u_{i}^{0}(x)\,,\quad u_{it}(T,x)=u_{i}^{1}(x)\,,
\quad x\in(0,\pi),\quad i=1,2\,.
\end{equation}
In the literature  coupled wave-wave equations were investigated
by studying  boundary stabilization, see  \cite{KR}.
The exact synchronization for a coupled system of wave equations
with Dirichlet boundary conditions was successfully treated by Li and Rao
\cite {LR}. They studied the $n-$dimensional case when the coupling matrix is very general. However, their method does not allow to get
precise estimates on the controllability time.

In \cite {Al} F. Alabau-Boussouira  considered a system where the coupling parameters are all equal,
obtaining  an observability inequality  for small coupling parameter  and large time $T$ and then, by duality,  an exact indirect controllability result.

In this paper we solve the  reachability problems for the coupled wave-wave 
 with an integro-differential term by the HUM method, see \cite{Lio1,Lio2,Lio3} and by means of non-harmonic analysis techniques. 
 In this framework  Ingham type estimates, see \cite{Ing}, play an important role. 
 We already used this approach to study the reachability for one equation, see \cite{LoretiSforza,LoretiSforza1} and to treat the case of a wave--Petrovsky system with a memory term, see \cite{LoretiSforza3}. For a different class of integral kernels see \cite{LPS} and for the hidden regularity in the case of general kernels see \cite{LoretiSforza4}.

 However the estimates obtained do not include the case wave-wave without memory as limit case
 as $\beta\to 0^+$
 \begin{equation}\label{eq:problem-usix}
\begin{cases}
\displaystyle 
u_{1tt} -u_{1xx}+au_2= 0
\\
\hskip4.5cm
\mbox{on}\quad
(0,T)\times (0,\pi),
\\
\displaystyle
u_{2tt} -u_{2xx}+bu_1= 0
\end{cases}
\end{equation} 
because, as  formulas \eqref {eq:lambda2} and \eqref {eq:lambda4} clearly show, the eigenvectors of the integro-differential operator
are not bounded as  $\beta\to 0^+ $.
  
  
  
The  method is based  on  a representation formula for the solution $(u_1,u_2)$, established
in Section 4
\begin{equation*}
\begin{split}
u_1(t) &=\sum_{n=1}^{\infty}\Big(C_{n}e^{i\omega_{n} t}+\overline{C_{n}}e^{-i\overline{\omega_{n}}t}
+R_{n}e^{r_{n} t}+D_{n}e^{i\zeta_{n} t}+\overline{D_{n}}e^{-i\overline{\zeta_{n}}t}\Big)
\,,
\\
u_2(t) &=\sum_{n=1}^{\infty}
\Big(d_nD_{n}e^{i\zeta_{n} t}+\overline{d_nD_{n}}e^{-i\overline{\zeta_{n}}t}
+c_nC_{n}e^{i\omega_{n} t}+\overline{c_nC_{n}}e^{-i\overline{\omega_{n}}t}\Big)
+\mathcal {E} e^{-\eta t}
\,,
\end{split}
\end{equation*}
where
\begin{equation*}
|\mathcal {E}|^2\le 
M \sum_{ n= 1}^{\infty}
\Big(|C_{n}|^2+ |d_nD_{n}|^2\Big),
\qquad
(M>0)
\,.
\end{equation*}
We will prove the following reachability result (see Theorem \ref{th:reachres}) where we will give an estimate of the control time.
\begin{theorem}
Let $\beta<1/2$. For  any $T>\frac{2\pi}{\sqrt{1-4\beta^2}}$ and
$
(u_{i}^{0},u_{i}^{1})\in  L^{2}(0,\pi)\times H^{-1}(0,\pi)
$,
$i=1,2$,
 there exist $g_i\in L^2(0,T)$, $i=1,2$, such that the weak solution  $(u_1,u_2)$ of system 
\begin{equation}\label{eq:problem-usix}
\begin{cases}
\displaystyle 
u_{1tt}(t,x) -u_{1xx}(t,x)+\beta\int_0^t\ e^{-\eta(t-s)} u_{1xx}(s,x)ds+au_2(t,x)= 0\,,
\\
\phantom{u_{1tt}(t,x) -u_{1xx}(t,x)+\int_0^t\ k(t-s) u_{1xx}(s,x)ds+}
t\in (0,T)\,,\,\,\, x\in(0,\pi)
\\
\displaystyle
u_{2tt}(t,x) -u_{2xx}(t,x)+bu_1(t,x)= 0
\,,
\end{cases}
\end{equation} 
with boundary conditions
\begin{equation}\label{eq:bound-u1r}
u_1(t,0)=u_2(t,0)=0\,,\quad u_1(t,\pi)=g_1(t)\,,\quad u_2(t,\pi)=g_2(t)\qquad t\in (0,T) 
 \,,
\end{equation}
and null initial values 
\begin{equation}
u_i(0,x)=u_{it}(0,x)=0\qquad  x\in(0,\pi)\,,\quad i=1,2,
\end{equation} 
verifies the final conditions
\begin{equation}\label{eq:findataT}
u_i(T,x)=u_{i}^{0}(x)\,,\quad u_{it}(T,x)=u_{i}^{1}(x)\,,
\quad x\in(0,\pi),
\qquad i=1,2\,.
\end{equation}
\end{theorem}
Due to the duality between controllability and observability we will first prove Ingham type inequalities (see Theorem \ref {th:inv.ingham1}).
\begin{theorem}\label{th:obsI}
Let $\{\om_n\}_{n\in\N}$, $\{r_n\}_{n\in\N}$  and
$\{\zeta_{n}\}_{n\in\N}$ be sequences of pairwise  distinct numbers
such that $\om_n\not= \zeta_m$, $\om_n\not=\overline{\zeta_m}$, 
$r_n\not= i\om_m$, $r_n\not= i\zeta_m$, $r_n\not=-\eta$, $\zeta_{n}\not=0$, for any $n\,,m\in\N$.
Assume that  there exist
$\gamma>0$, $\alpha,\chi\in\R$, $n'\in\N$, $\mu>0$, $\nu> 1/2$, 
such that
\begin{equation*}
\liminf_{n\to\infty}({\Re}\om_{n+1}-{\Re}\om_{n})=\liminf_{n\to\infty}({\Re}\zeta_{n+1}-{\Re} \zeta_{n})=\gamma\,,
\end{equation*}
\begin{equation*}
\begin{split}
\lim_{n\to\infty}{\Im}\om_n&=\alpha>0
\,,
\\
\lim_{n\to\infty}r_n&=\chi<0\,,
\\
\lim_{n\to\infty}\Im \zeta_{n}&=0\,,
\end{split}
\end{equation*}
\begin{equation*}
|d_n|\asymp|\zeta_n|
\,,
\qquad
|c_n|\le\frac{M}{|\omega_n|}\,,
\end{equation*}
\begin{equation*}
|R_n|\le \frac{\mu}{n^{\nu}}\Big(|C_{n}|^2+ |d_nD_{n}|^2\Big)^{1/2}\,\quad\forall\ n\ge n'\,,
\qquad
|R_n|\le \mu\Big(|C_{n}|^2+ |d_nD_{n}|^2\Big)^{1/2}\,\quad\forall\ n\le n'\,.
\end{equation*}
Then, for $\gamma>4\alpha$ and
$T>\frac{2\pi}{\sqrt{\gamma^2-16\alpha^2}}$ we have 
 \begin{equation}\label{eq:inv.ingham}
 \int_{0}^{T} \big(|u_1(t)|^2+|u_2(t)|^2\big)\ dt
 \asymp
\sum_{ n= 1}^{\infty}
\Big(|C_{n}|^2+ |d_nD_{n}|^2\Big) \,.
\end{equation}
\end{theorem}
The observability time may be improved making an extra assumption
on the initial data. Indeed, if we assume the condition $|C_n|\le M |d_nD_{n}|$ on the coefficients of the series
 instead of  $\gamma>4\alpha$, then we can make use of Theorem \ref{th:extracoe} instead of Theorem \ref{th:gamma>4alpha}, obtaining the observability estimates for $T>\frac{2\pi}{\gamma}$ (see Theorem \ref{th:inv.ingham11}).

\begin{theorem}
Let  assume the hypotheses of Theorem \ref{th:obsI} and the condition 
\begin{equation}
|C_n|\le M |d_nD_{n}|
\,.
\end{equation}
Then, for 
$T>\frac{2\pi}{\gamma}$ we have 
 \begin{equation}\label{eq:inv.ingham11}
 \int_{0}^{T} \big(|u_1(t)|^2+|u_2(t)|^2\big)\ dt
 \asymp
\sum_{ n= 1}^{\infty}
\Big(|C_{n}|^2+ |d_nD_{n}|^2\Big) \,.
\end{equation}

\end{theorem}
The plan of our paper is the following. In Section 2 we  give some preliminary results. In Section 3 we describe the Hilbert Uniqueness Method. In Section 4 we carry out a detailed spectral analysis to give a representation formula for the solution of the wave-wave coupled system with memory.
In Section 5 we prove the observability estimates.
Finally, in Section 6 we give a reachability result for the coupled system with memory.

\section{Preliminaries}
Throughout the paper, we will adopt the convention to write $F\asymp G$ if there exist two positive constants $c_1$ and $c_2$ such that
$c_1F\le G\le c_2F$.

Let $X$ be a real Hilbert space with scalar product
$\langle \cdot \, ,\, \cdot \rangle$ and norm $\| \cdot \|$. 
For any
$T\in\, (0,
\infty]$   we denote by
$L^1(0,T;X)$ the usual spaces of measurable
functions
$v:(0,T)\to X$ such that one has

$$
\|v\|_{1,T}:=\int_0^T \|v(t)\|\,dt<\infty\,.
$$
 We shall use the shorter notation $\|v\|_1$ for
$\|v\|_{1,\infty}$.
We denote by $L_{loc}^1 (0,\infty;X)$ the space of functions
belonging to
$L^1(0,T;X)$ for any $T\in (0,\infty)$.
In the case of $X=\R$, we will use the abbreviations
$L^1(0,T)$ and
$L_{loc}^1(0,\infty)$ to denote the spaces $L^1(0,T;\R)$
and
$L_{loc}^1(0,\infty;\R)$,
respectively.

Classical results for integral equations
(see, e.g., \cite[Theorem 2.3.5]{GLS})
ensure that, for any  kernel $k\in L_{loc}^1(0,\infty)$ and 
$\psi\in L_{loc}^1 (0,\infty;X)$,  the problem
\begin{equation}\label{integral}
\varphi(t)-k*\varphi(t)=\psi(t),\qquad t\ge 0\,,
\end{equation}
admits a unique solution $\varphi\in L_{loc}^1(0,\infty;X)$. In particular, if we take $\psi=k$ in \eqref{integral}, 
we can consider the unique
solution $\varrho_k\in L_{loc}^1(0,\infty)$ of 
\begin{equation*}
\varrho_k (t)-k*\varrho_k (t)=k (t),\qquad t\ge 0\,.
\end{equation*}
Such a solution is called the {\em resolvent kernel} of $k$.
Furthermore, for any $\psi$ the solution $\varphi$ of (\ref{integral}) is given by the
variation of constants formula
\begin{equation*}
\varphi(t)=\psi(t)+\varrho_k *\psi(t),\qquad t\ge 0\,,
\end{equation*}
where $\varrho_k$ is the resolvent kernel of $k$.

We recall some results concerning integral equations in case of decreasing exponential kernels, see for example \cite[Corollary 2.2]{LoretiSforza1}.
\begin{proposition}\label{pr:unicita}
For $0<\beta<\eta$ and $T>0$ the following properties hold true.
\begin{description}
\item [(i)] The resolvent kernel of $k(t)=\beta e^{-\eta t}$ is $\varrho_k(t)=\beta e^{(\beta-\eta)t}$.

\item [(ii)]Given $\psi\in
L_{loc}^1 (-\infty,T;X)$, a function $\varphi\in L_{loc}^1 (-\infty,T;X)$ is a solution of
\begin{eqnarray*}
\varphi(t)-\beta\int_t^T\ e^{-\eta(s-t)}\varphi(s)ds=\psi(t)
\qquad t\le T\,,
\end{eqnarray*}
if and only if
\begin{eqnarray*}
\varphi(t)=\psi(t)+\beta\int_{t}^Te^{(\beta-\eta)(s-t)}\psi(s)\ ds
\qquad t\le T\,.
\end{eqnarray*}
Moreover, there exist two positive constants $c_1\,,c_2$ depending on $\beta,\eta,T$ such that
\begin{equation}\label{eq:unicita}
c_1\int_0^T|\varphi(t)|^2\ dt
\le\int_0^T|\psi(t)|^2\ dt
\le c_2\int_0^T|\varphi(t)|^2\ dt\,.
\end{equation}

\end{description}

\end{proposition}

We state and prove a result, that will allow us  to give an equivalent way to write the solution of our problem.

\begin{lemma}\label{le:fifth}
Given $\la\,,\beta\,,\eta\in\R$, $a\in\R\setminus\{0\}$ and $b\in\R$, a couple $(f, g)$ of scalar functions defined on the interval $[0,\infty)$ is a solution of the system
\begin{equation}\label{eq:system0}
\begin{cases}\displaystyle
f{''}+\la f-\la\be \int_0^t e^{-\eta(t-s)}f(s) ds+ag=0\,,
\\
\hskip8cm
t\ge 0,
\\\displaystyle
g{''}+\la g+bf=0
\,,
\end{cases}
\end{equation}
if and only if $f$ is  a solution of the equation
\begin{equation}\label{eq:fifth}
\displaystyle 
f^{(5)}+\eta f^{(4)}+2\la f{'''}+\la (2\eta- \be) f''+(\la^2-ab)f'
+(\la^2 (\eta- \be)-\eta ab)f=0,\quad t\ge 0,
\end{equation}
the condition 
\begin{equation}\label{eq:fifth1}
\displaystyle
 f^{(4)}(0)=-2\la f{''}(0)+\la \be f'(0)+(ab-\eta\la \be -\la^2)f(0)
\end{equation}
is satisfied and $g$ is given by
\begin{equation}\label{eq:Ag}
g=-\frac1a\Big(f{''}+\la f-\la\be \int_0^t e^{-\eta(t-s)}f(s) ds\Big)\,.
\end{equation}
\end{lemma}
\begin{Proof}
Let $(f, g)$ be a solution of (\ref{eq:system0}). 
Differentiating the first equation in (\ref{eq:system0}), we get
\begin{equation}\label{eq:0gprime}
f{'''}+\la f{'}
+\eta\la\be \int_0^t e^{-\eta(t-s)}f(s) ds-\la \be f+ag'=0\,,
\end{equation}
whence
\begin{equation}\label{eq:gprime0}
ag'(0)=-f{'''}(0)-\la f{'}(0)+\la \be f(0)\,.
\end{equation}
Substituting in \eqref{eq:0gprime} the  identity
\begin{equation*}
\la\be \int_0^t e^{-\eta(t-s)}f(s) ds=f{''}+\la f+ag\,,
\end{equation*}
we obtain
\begin{equation}\label{eq:gprime}
f{'''}+\eta f{''}+\la f{'}+\la (\eta- \be) f+ag'+\eta ag=0\,.
\end{equation}
Differentiating yet again, we have
\begin{equation*}
f^{(4)}+\eta f{'''}+\la f{''}+\la (\eta- \be) f'+ag''+\eta ag'=0\,,
\end{equation*}
whence, by using the second equation in (\ref{eq:system0}), 
that is
$ag{''}=-abf-\la ag$,
we get
\begin{equation}\label{eq:fourth}
f^{(4)}+\eta f{'''}+\la f{''}+\la (\eta- \be) f'-abf+\eta ag'-\la ag=0\,.
\end{equation}
Thanks to \eqref{eq:gprime0} and $ag(0)=-f{''}(0)-\la f(0)$, we have
\begin{multline*}
f^{(4)}(0)=-\eta f{'''}(0)-\la f{''}(0)-\la (\eta- \be) f'(0)+abf(0)-\eta ag'(0)+\la ag(0)\\
=-\eta f{'''}(0)-\la f{''}(0)-\la (\eta- \be) f'(0)+abf(0)+\eta f{'''}(0)\\
+\eta\la f{'}(0)-\eta\la \be f(0)-\la f{''}(0)-\la^2 f(0)
\\
=-2\la f{''}(0)+\la \be f'(0)+(ab-\eta\la \be -\la^2)f(0)
\,,
\end{multline*}
so formula (\ref{eq:fifth1}) for $f^{(4)}(0)$ holds true.
Moreover,
by differentiating  \eqref{eq:fourth} we obtain
\begin{equation*}
f^{(5)}+\eta f^{(4)}+\la f{'''}+\la (\eta- \be) f''-abf'+\eta ag''-\la ag'=0\,.
\end{equation*}
By using again
$g{''}=-bf-\la g$
we get 
\begin{equation*}
f^{(5)}+\eta f^{(4)}+\la f{'''}+\la (\eta- \be) f''-abf'-\eta abf-\la ag'-\eta\la  ag=0\,.
\end{equation*}
From \eqref{eq:gprime} it follows 
\begin{equation*}
-ag'-\eta ag=f{'''}+\eta f{''}+\la f{'}+\la (\eta- \be) f\,,
\end{equation*}
and hence we have
\begin{equation*}
f^{(5)}+\eta f^{(4)}+2\la f{'''}+\la (2\eta- \be) f''
+(\la^2-ab)f'+(\la^2 (\eta- \be)-\eta ab)f
=0\,,
\end{equation*}
that is
$f$ is a solution of the differential equation (\ref{eq:fifth}). 
Finally,
from the first equation in (\ref{eq:system0}) we deduce that
$g$ is given by \eqref{eq:Ag}.

Conversely, if $f$ satisfies $(\ref{eq:fifth})-(\ref{eq:fifth1})$,  
multiplying the differential equation
 by $e^{\eta t}$ and integrating from $0$ to $t$, we obtain 
\begin{multline*}
\int_0^t e^{\eta s}f^{(5)}(s)\ ds+
\eta \int_0^t e^{\eta s}f^{(4)}(s)\ ds
+2\la\int_0^t e^{\eta s}f{'''}(s)\ ds
 +2\eta\la\int_0^t e^{\eta s}f{''}(s)\ ds 
 \\
-\la\be\int_0^t e^{\eta s}f{''}(s)\ ds
+(\la^2-ab)\int_0^t e^{\eta s}f{'}(s)\ ds
+(\la^2(\eta-\be)-\eta ab)\int_0^t e^{\eta s}f(s)\ ds=0\,.
\end{multline*}
Integrating by parts the first, the third, the fifth and the sixth integral, we have
\begin{multline*}
e^{\eta t}f^{(4)}-f^{(4)}(0)
+2\la e^{\eta t}f{''} -2\la f{''}(0)
-\la\be e^{\eta t}f{'}+\la\be f{'}(0)+\eta\la\be e^{\eta t}f{} \\
-\eta\la\be  f{}(0)-\eta^2\la\be\int_0^t e^{\eta s}f{}(s)\ ds
+(\la^2-ab)e^{\eta t}f-(\la^2-ab)f(0)
-\la^2\be\int_0^t e^{\eta s}f(s)\ ds=0\,.
\end{multline*}
Using the condition (\ref{eq:fifth1}) and multiplying by $e^{-\eta t}$, we obtain
\begin{multline}\label{eq:fourthbis}
f^{(4)}
+2\la f{''} 
-\la\be f{'}+\eta\la\be f{} 
-\eta^2\la\be\int_0^t e^{-\eta(t- s)}f{}(s)\ ds\\
+(\la^2-ab)f
-\la^2\be\int_0^t e^{-\eta(t- s)}f(s)\ ds=0\,.
\end{multline}
Moreover, by \eqref{eq:Ag} it follows
\begin{equation*}
ag'=-f{'''}-\la f'+\la\beta f-\eta\la\be \int_0^t e^{-\eta(t-s)}f(s) ds\,,
\end{equation*}
and hence
\begin{equation*}
ag''=-f^{(4)}-\la f''+\la\beta f'-\eta\la\be f +\eta^2\la\be\int_0^t e^{-\eta(t-s)}f(s) ds\,.
\end{equation*}
Therefore, thanks to the previous identity and \eqref{eq:fourthbis} we have
\begin{equation*}
ag''=
\la f{''} 
+(\la^2-ab)f
-\la^2\be\int_0^t e^{-\eta(t- s)}f(s)\ ds\,,
\end{equation*}
whence, in view of \eqref{eq:Ag} we get
\begin{equation*}
ag''=-\la ag-abf\,.
\end{equation*}
Finally, by \eqref{eq:Ag} and the above equation, it follows that the couple $(f, g)$ is a solution of the system \eqref{eq:system0}.
\end{Proof}

The following lemma is analogous to that of \cite[Lemma 2.3]{LoretiSforza1}. For the reader's convenience we prefer to state and prove it the same.
\begin{lemma}\label{le:third}
Given $\la\,,\beta\,,\eta\in\R$ and $h\in C(\R)$, if  $g\in C^3(\R)$ 
is  a solution of the third order differential equation
\begin{equation}\label{third}
g{'''}+
\eta g{''}+\la g{'}+\la (\eta-\be)g=h\,\qquad \mbox{in}\,\,\,\R
\,,
\end{equation}
then $g$ is also a solution of the integro-differential equation 
\begin{equation}\label{eq:second}
g{''}
+\la g-\la\be \int_0^t e^{-\eta(t-s)}g(s) ds=e^{-\eta t}(g{''}(0)+\la g(0))+ \int_0^t e^{-\eta(t-s)}h(s) ds\,\qquad t\in\R\,.
\end{equation}

\end{lemma}
\begin{Proof}
Multiplying the differential equation $(\ref{third})$ by $e^{\eta t}$ and integrating from $0$ to $t$, we obtain 
\begin{equation*}
\int_0^t e^{\eta s}g{'''}(s)\ ds+
\eta \int_0^t e^{\eta s}g{''}(s)\ ds+\la\int_0^t e^{\eta s}g{'}(s)\ ds  
+\la (\eta-\be)\int_0^t e^{\eta s}g(s)\ ds=\int_0^t e^{\eta s}h(s)\ ds\,.
\end{equation*}
Integrating by parts the first term and the third one, we have
\begin{equation*}\label{third1} 
e^{\eta t}g{''}-g{''}(0)+\la e^{\eta t} g-\la g(0)-
\la\be \int_0^t e^{\eta s}g(s)\ ds=\int_0^t e^{\eta s}h(s)\ ds\,.
\end{equation*}
Finally, if we multiply by $e^{-\eta t}$, then we obtain  $(\ref{eq:second})$. 
\end{Proof}

\section{The Hilbert Uniqueness Method}\label{se:HUM}
For reader's convenience,  in this section we will describe the Hilbert Uniqueness Method for coupled wave equations with a memory term.
For another  approach based
on the ontoness  of the solution operator, see e.g. \cite{LasT, T1}.

Given $k\in L_{loc}^1(0,\infty)$ and $a\,,b\in\R$,
we consider the following coupled system:
\begin{equation}\label{eq:problem-u}
\begin{cases}
\displaystyle 
u_{1tt}(t,x) -u_{1xx}(t,x)+\int_0^t\ k(t-s) u_{1xx}(s,x)ds+au_2(t,x)= 0\,,
\\
\phantom{u_{1tt}(t,x) -u_{1xx}(t,x)+\int_0^t\ k(t-s) u_{1xx}(s,x)ds+au_2(t,x)= 0\,,\qquad}
t\in (0,T)\,,\quad x\in(0,\pi)
\\
\displaystyle
u_{2tt}(t,x) -u_{2xx}(t,x)+bu_1(t,x)= 0
\,,
\end{cases}
\end{equation}
subject to the boundary conditions
\begin{equation}\label{eq:bound-u2}
u_1(t,0)=u_2(t,0)=0\,,\quad u_1(t,\pi)=g_1(t)\,,\quad u_2(t,\pi)=g_2(t)\qquad t\in (0,T) \,,
\end{equation}
and with null initial conditions 
\begin{equation}\label{eq:initc}
u_i(0,x)=u_{it}(0,x)=0\qquad  x\in(0,\pi),\quad i=1,2\,.
\end{equation} 

For a reachability problem we mean the following: given $T>0$ and taking
$
(u_{i}^{0},u_{i}^{1})
$, $i=1,2$,
in a suitable space, that we will introduce later,
find $g_i\in L^2(0,T)$, $i=1,2$ such that the weak
solution $u$ of problem  \eqref{eq:problem-u}-\eqref{eq:initc}
satisfies the final conditions
\begin{equation}\label{eq:problem-u1}
u_i(T,x)=u_{i}^{0}(x)\,,\quad u_{it}(T,x)=u_{i}^{1}(x)\,,
\quad x\in(0,\pi),\quad i=1,2\,.
\end{equation}
One can solve such reachability problems by the HUM method. To see that, we proceed as follows.

Given
$(z_{i}^{0},z_{i}^{1})\in (C^\infty_c(0,\pi))^2$, $i=1,2$,
we introduce the {\it adjoint} system of (\ref{eq:problem-u}), that is 
\begin{equation}\label{eq:adjoint}
\begin{cases}
\displaystyle 
z_{1tt}(t,x) -z_{1xx}(t,x)+\int_t^T\ k(s-t) z_{1xx}(s,x)ds+bz_2(t,x)= 0\,,\\
\hskip9.5cm
t\in (0,T)\,,\quad x\in(0,\pi)
\\
\displaystyle
z_{2tt}(t,x) -z_{2xx}(t,x)+az_1(t,x)= 0\,,
\\
z_i(t,0)=z_i(t,\pi)=0\qquad t\in [0,T], \quad i=1,2,
\end{cases}
\end{equation}
with  final data 
\begin{equation} \label{eq:final}
z_i(T,\cdot)=z_{i}^{0}\,,\quad z_{it}(T,\cdot)=z_{i}^{1}\,,\quad i=1,2
\,.
\end{equation} 
The above problem is well-posed, see e.g. \cite{Pruss}. 
Thanks to the regularity of the final data, the solution  $(z_1,z_2)$ of \eqref{eq:adjoint}--\eqref{eq:final}
is regular enough to 
consider  the  nonhomogeneous problem
\begin{equation}\label{eq:phi}
\left \{\begin{array}{l}\displaystyle
\phi_{1tt}(t,x) -\phi_{1xx}(t,x)+\int_0^t\ k(t-s) \phi_{1xx}(s,x)ds+a \phi_ 2(t,x)= 0
\\
\hskip9.5cm  t\in (0,T)\,,\quad x\in(0,\pi)\,,
\\
\displaystyle
\phi_{2tt}(t,x) - \phi_{2xx}(t,x)+b \phi_ 1(t,x)= 0
\\
\\
\phi_i(0,x)= \phi_{it}(0,x)=0\qquad  x\in(0,\pi)\,, \quad i=1,2,
\\
\\
\displaystyle
\phi_1(t,0)=0\,,\quad
\phi_1(t,\pi)=z_{1x}(t,\pi)-\int_t^T\ k(s-t)z_{1x}(s,\pi)ds
\\
\hskip9.5cm  t\in [0,T],
\\
\phi_2(t,0)=0\,,\quad 
\phi_2(t,\pi)=z_{2x}(t,\pi).
\end{array}\right .
\end{equation} 
As in the non-integral case, it can be proved  that  problem \eqref{eq:phi} admits a unique solution
$(\phi_1,\phi_2)$.
So, we can introduce the following linear operator: for any  $(z_{i}^{0},z_{i}^{1})\in \big(C^\infty_c(0,\pi)\big)^2$, $i=1,2$,
we define
\begin{equation}\label{eq:psi0}
\Psi(z_{1}^{0},z_{1}^{1},z_{2}^{0},z_{2}^{1})=(-\phi_{1t}(T,\cdot),\phi_{1}(T,\cdot),-\phi_{2t}(T,\cdot),\phi_{2}(T,\cdot))
\,.
\end{equation}
For any $(\xi_{i}^{0},\xi_{i}^{1})\in \big(C^\infty_c(0,\pi)\big)^2$, $i=1,2$, let $(\xi_1,\xi_2)$ be the solution of 
\begin{equation}\label{eq:adjoint10}
\left \{\begin{array}{l}\displaystyle
\xi_{1tt}(t,x) -\xi_{1xx}(t,x)+\int_t^T\ k(s-t) \xi_{1xx}(s,x)ds+b\xi_2(t,x)= 0
\\
\hskip9.5cm
t\in (0,T),\quad x\in(0,\pi),
\\
\displaystyle
\xi_{2tt}(t,x) -\xi_{2xx}(t,x)+a\xi_1(t,x)= 0
\\
\\
\xi_i(t,0)=\xi_i(t,\pi)=0
\qquad t\in [0,T],
\\
\hskip6cm
\quad i=1,2,
\\
\xi_i(T,\cdot)=\xi_{i}^{0}\,,\quad \xi_{it}(T,\cdot)=\xi_{i}^{1}\,.
\end{array}\right .
\end{equation} 
We will prove that
\begin{multline}\label{eq:psi}
\langle\Psi(z_{1}^{0},z_{1}^{1},z_{2}^{0},z_{2}^{1}),(\xi_{1}^{0},\xi_{1}^{1},\xi_{2}^{0},\xi_{2}^{1})\rangle_{L^2}
\\
=\int_0^T\phi_1(t,\pi)\Big(\xi_{1x}(t,\pi)-\int_t^T\ k(s-t)\ \xi_{1x}(s,\pi)\ ds\Big) \ dt
+\int_0^T\phi_{2}(t,\pi)\xi_{2x}(t,\pi)\ dt
 \,. 
\end{multline}
To this end, we multiply the first equation in (\ref{eq:phi}) by $\xi_1$ and integrate 
on $[0,T]\times[0,\pi]$, so we have
\begin{multline*}
\int_0^\pi \int_0^T\phi_{1tt}(t,x)\xi_1(t,x)\ dt \ dx
-\int_0^T\int_0^\pi\phi_{1xx}(t,x)\xi_1(t,x)\ dx\ dt
\\
+\int_0^\pi\int_0^T\int_0^t\ k(t-s)\phi_{1xx}(s,x)\ ds\ \xi_1(t,x)\ dt\ dx
+a\int_0^T\int_0^\pi \phi_{2}(t,x)\xi_1(t,x)\ dx \ dt=0\,. 
\end{multline*}
If we take into account that 
\begin{equation*}
\int_0^T\int_0^t\ k(t-s)\phi_{1xx}(s,x)\ ds\ \xi_1(t,x)\ dt  =
\int_0^T\phi_{1xx}(s,x)\int_s^T\ k(t-s)\ \xi_1(t,x)\ dt \ ds 
\end{equation*}
and
integrate by parts,  then we have
\begin{multline*}
 \int_0^\pi\big(\phi_{1t}(T,x)\xi_{1}^{0}(x)- \phi_1(T,x)\xi_{1}^{1}(x)\big)\ dx
 +\int_0^\pi \int_0^T\phi_1(t,x)\xi_{1tt}(t,x)\ dt \ dx
 \\
+\int_0^T\phi_1(t,\pi)\xi_{1x}(t,\pi)\ dt -\int_0^T\int_0^\pi\phi_1(t,x)\xi_{1xx}(t,x)\ dx\ dt
 \\
-\int_0^T\phi_{1}(s,\pi)\int_s^T\ k(t-s)\ \xi_{1x}(t,\pi)\ dt \ ds 
+\int_0^\pi\int_0^T\phi_1(s,x)\int_s^T\ k(t-s)\ \xi_{1xx}(t,x)\ dt \ ds \ dx
\\
+a\int_0^T\int_0^\pi \phi_{2}(t,x)\xi_1(t,x)\ dx \ dt=0\,. 
\end{multline*}
As a consequence of the above equation and
\begin{equation*}
\xi_{1tt} -\xi_{1xx}+\int_t^T\ k(s-t) \xi_{1xx}(s,\cdot)ds=-b\xi_2\,,
\end{equation*}
we obtain
\begin{multline}\label{eq:xi1}
 \int_0^\pi\big(\phi_{1t}(T,x)\xi_{1}^{0}(x)- \phi_1(T,x)\xi_{1}^{1}(x)\big)\ dx
+\int_0^T\phi_1(t,\pi)\Big(\xi_{1x}(t,\pi)-\int_t^T\ k(s-t)\ \xi_{1x}(s,\pi)\ ds\Big) \ dt 
\\
+\int_0^T\int_0^\pi \big(a\phi_{2}(t,x)\xi_1(t,x)-b\phi_1(t,x)\xi_{2}(t,x)\big)\ dx \ dt=0\,. 
\end{multline}
In a similar way, we multiply the second equation in (\ref{eq:phi}) by $\xi_2$ and integrate by parts
on $[0,T]\times[0,\pi]$ to get
\begin{multline*}
 \int_0^\pi\big(\phi_{2t}(T,x)\xi_{2}^{0}(x)- \phi_2(T,x)\xi_{2}^{1}(x)\big)\ dx
 +\int_0^\pi \int_0^T\phi_2(t,x)\xi_{2tt}(t,x)\ dt \ dx
 \\
+\int_0^T\phi_{2}(t,\pi)\xi_{2x}(t,\pi)\ dt 
-\int_0^T\int_0^\pi\phi_2(t,x)\xi_{2xx}(t,x)\ dx\ dt
+b\int_0^T\int_0^\pi \phi_{1}(t,x)\xi_2(t,x)\ dx \ dt=0\,, 
\end{multline*}
whence, in virtue of
\begin{equation*}
\xi_{2tt} -\xi_{2xx}=-a\xi_1\,,
\end{equation*}
we get
\begin{multline}\label{eq:xi2}
 \int_0^\pi\big(\phi_{2t}(T,x)\xi_{2}^{0}(x)- \phi_2(T,x)\xi_{2}^{1}(x)\big)\ dx
+\int_0^T\phi_{2}(t,\pi)\xi_{2x}(t,\pi)\ dt 
\\
+\int_0^T\int_0^\pi\big(b \phi_{1}(t,x)\xi_2(t,x)-a\phi_{2}(t,x)\xi_1(t,x)\big)\ dx \ dt=0\,. 
\end{multline}
If we sum equations \eqref{eq:xi1} and \eqref{eq:xi2}, then we have
\begin{multline}\label{eq:prenorm}
\langle\Psi(z_{1}^{0},z_{1}^{1},z_{2}^{0},z_{2}^{1}),(\xi_{1}^{0},\xi_{1}^{1},\xi_{2}^{0},\xi_{2}^{1})\rangle_{L^2}
\\
= \int_0^\pi\big(-\phi_{1t}(T,x)\xi_{1}^{0}(x)+ \phi_1(T,x)\xi_{1}^{1}(x)-\phi_{2t}(T,x)\xi_{1}^{0}(x)+ \phi_2(T,x)\xi_{1}^{1}(x)\big)\ dx
 \\
 =\int_0^T\phi_1(t,\pi)\Big(\xi_{1x}(t,\pi)-\int_t^T\ k(s-t)\ \xi_{1x}(s,\pi)\ ds\Big) \ dt
 +\int_0^T\phi_{2}(t,\pi)\xi_{2x}(t,\pi)\ dt \,,
 \end{multline}
that is, \eqref{eq:psi} holds true.

Taking $\xi_{i}^{0}=z_{i}^{0}$ and $\xi_{i}^{1}=z_{i}^{1}$, $i=1,2$, in  (\ref{eq:psi}) yields
\begin{multline}\label{eq:psi1}
\langle\Psi(z_{1}^{0},z_{1}^{1},z_{2}^{0},z_{2}^{1}),(z_{1}^{0},z_{1}^{1},z_{2}^{0},z_{2}^{1})\rangle_{L^2}
\\
=
\int_0^T\Big|z_{1x}(t,\pi)-\int_t^T\ k(s-t)\ z_{1x}(s,\pi)\ ds\Big|^2 \ dt
 +\int_0^T\big|z_{2x}(t,\pi)\big|^2\ dt\,. 
\end{multline}
As a consequence, we can introduce a semi-norm on the space $\big(C^\infty_c(\Omega)\big)^4$. Indeed,
 for $(z_{i}^{0},z_{i}^{1})\in \big(C^\infty_c(\Omega)\big)^2$, $i=1,2$, we define
\begin{multline}\label{eq:normF}
\|(z_{1}^{0},z_{1}^{1},z_{2}^{0},z_{2}^{1})\|_{F}:=
\displaystyle\Big(
\int_0^T\Big|z_{1x}(t,\pi)-\int_t^T\ k(s-t)\ z_{1x}(s,\pi)\ ds\Big|^2 \ dt
 +\int_0^T\big|z_{2x}(t,\pi)\big|^2\ dt
\Big)^{1/2}.
\end{multline}
In view of 
Proposition \ref{pr:unicita},  $\|\cdot\|_{F}$ is a norm if and only if the following uniqueness theorem holds.
\begin{theorem}\label{th:uniqueness}
If $(z_1,z_2)$ is the solution of problem {\rm (\ref{eq:adjoint})--(\ref{eq:final})} such that
$$
z_{1x}(t,\pi)=z_{2x}(t,\pi)=0\,,\qquad \forall t\in [0,T]\,,
$$
then 
$$
z_1(t,x)=z_2(t,x)= 0 \qquad\forall (t,x)\in [0,T]\times[0,\pi]\,.
$$
\end{theorem}
If we are able to establish Theorem \ref{th:uniqueness}, then we can define the Hilbert space $F$ as the completion of $ \big(C^\infty_c(\Omega)\big)^4$ for
the norm (\ref{eq:normF}). Moreover, the operator $\Psi$ extends uniquely to a continuous operator, denoted again by $\Psi$, from $F$ to the dual
space $F'$ in such a way that   $\Psi:F\to F'$ is an isomorphism.

In conclusion, if we prove Theorem \ref{th:uniqueness} and,  for example,
$F=\big(H^1_0(0,\pi)\times L^2(0,\pi)\big)^2$
with the equivalence of the respective norms, then, taking $(u_{i}^{0},u_{i}^{1})\in  L^{2}(0,\pi)\times H^{-1}(0,\pi)$, $i=1,2$, we can solve the reachability problem
\eqref{eq:problem-u}--\eqref{eq:problem-u1}.

\section{Representation of the solution as Fourier series}\label{se:}
\subsection{Spectral analysis}\label{se:specan}

The aim of this section will be to give a complete spectral analysis for the coupled system. 

We will recast our system of coupled wave equations with a memory term in an abstract setting. Indeed, we consider
a self-adjoint positive linear
 operator 
$L:D(L)\subset H\to H$ on a Hilbert space $H$ with dense domain $D(L)$. We denote by $\{\la_n\}_{n\ge1}$ a strictly increasing sequence  of  eigenvalues for the operator $L$ with
$\la_n>0$ and $\la_n\to\infty$ and we assume that the sequence of the corresponding eigenvectors $\{w_n\}_{n\ge1}$ constitutes a Hilbert basis for $H$.

We fix two real numbers $a\not=0$, $b$
and consider the following coupled system:
\begin{equation}\label{eq:system}
\begin{cases}
\displaystyle 
u_1''(t) +Lu_1(t)-\beta\int_0^t\ e^{-\eta(t-s)}L u_1(s)ds+au_2(t)= 0
\\
\hskip9cm  t\ge 0,
\\
\displaystyle
u_2''(t) +Lu_2(t)+bu_1(t)= 0
\\
u_i(0)=u_{i}^{0}\,,\quad u'_i(0)=u_{i}^{1}\,,
\quad i=1,2
\,.
\end{cases}
\end{equation}
If we take the initial data $(u_{i}^{0},u_{i}^{1})$, $i=1,2$, belonging to $D(\sqrt{L})\times H$, then we can expand them according to the eigenvectors $w_n$ to obtain:
\begin{equation}\label{eq:v0}
\begin{split}
& u_{i}^{0}=\sum_{n=1}^{\infty}\alpha_{in}w_{n}\,,\qquad\quad\alpha_{in}=
\langle u_{i}^{0},w_n\rangle \,,
\quad
\|u_i^0\|^2_{D(\sqrt{L})}:=\sum_{n=1}^{\infty}\alpha_{in}^2\lambda_n\,,
\\
& u_{i}^{1}=\sum_{n=1}^{\infty}\rho_{in}w_{n}\,,\qquad\quad\rho_{in}=\langle u_{i}^{1},w_n\rangle \,,\quad
\|u_i^1\|^2_{H}:=\sum_{n=1}^{\infty}\rho_{in}^2\,.
\end{split} 
\end{equation} 
Our target is to write the components $u_{1}, u_{2}$ of the solution of system \eqref{eq:system} as sums of series, that is
\begin{equation*}
u_{i}(t)=\sum_{n=1}^{\infty}f_{in}(t)w_{n}\,,
\qquad
f_{in}(t)=\langle u_{i}(t),w_n\rangle\,,\quad i=1,2
\,.
\end{equation*}
To this end, we put the above  expressions for $u_{1}$ and $u_{2}$
into \eqref{eq:system} and multiply  by $w_n$,  
so  for any $n\in\N$ $(f_{1n},f_{2n})$ is the solution of the system
\begin{equation}\label{eq:secondsys}
\begin{cases}\displaystyle
f_{1n}''
+\la_{n}f_{1n}-\beta\la_{n}\int_0^t e^{-\eta(t-s)}f_{1n}(s) ds
+af_{2n}=0,
\\\displaystyle
f_{2n}''+\la_n f_{2n}+bf_{1n}=0
\,,
\\
f_{in}(0)=\a_{in}\,, \quad f_{in}'(0)=\rho_{in}\,,
\quad i=1,2
\,.
\end{cases}
\end{equation}
Thanks to lemma \ref{le:fifth} with $\la=\la_n$, $(f_{1n},f_{2n})$ is the solution of problem \eqref{eq:secondsys} if and only if  $f_{1n}$ is the solution of the Cauchy problem
\begin{equation}\label{eq:third}
\begin{cases}
\displaystyle 
f_{1n}^{(5)}+\eta f_{1n}^{(4)}+2\la_{n} f_{1n}'''+\la_{n} (2\eta-\be)f_{1n}''+(\la_{n}^2-ab)f_{1n}'
+(\la_{n}^2(\eta-\be)-\eta ab)f_{1n}=0
\qquad t\ge 0\,,
\\
 f_{1n}(0)=\alpha_{1n},
\\
 f_{1n}'(0)=\rho_{1n},
 \\
 f_{1n}''(0)=-\la_{n}\alpha_{1n}-a\alpha_{2n},
 \\
 f_{1n}'''(0)=-\la_{n} \rho_{1n}+ \beta\la_{n}\alpha_{1n}-a\rho_{2n},
 \\
 f_{1n}^{(4)}(0)
 =(\la_{n}^2-\eta \be\la_n+ab)\alpha_{1n}+2a\la_{n}\alpha_{2n}+\be\la_{n} \rho_{1n}\,,
\end{cases}
\end{equation} 
and $f_{2n}$ is given by
\begin{equation*}
f_{2n}=-\frac1a\Big(f_{1n}''+\la_n f_{1n}-\beta\la_{n} \int_0^t e^{-\eta(t-s)}f_{1n}(s) ds\Big)\,.
\end{equation*}
If we introduce the linear operator 
$\Upsilon_n$ defined by
\begin{equation}\label{eq:f2j'}
\Upsilon_n(v)(t):=
-\frac1a\Big(v''(t)
+\la_{n}v(t)-\beta\la_{n}\int_0^t e^{-\eta(t-s)}v(s) ds\Big)
\qquad t\ge0
\,,
\end{equation} 
then $f_{2n}$ can be written as
\begin{equation}\label{eq:f2j}
f_{2n}(t)
=\Upsilon_n(f_{1n})(t)
\qquad t\ge0
\,.
\end{equation}
We also note that for any $z\in\C$
\begin{equation}\label{eq:Yexp}
\Upsilon_n(e^{zt})=
-\frac1a
\Big[\Big(z^2+\la_n-\frac{\beta\la_n}{\eta+z}\Big)e^{z t}
+\frac{\beta \la_n}{\eta+z}e^{-\eta t}
\Big]\,.
\end{equation}
\subsection{The fifth order ordinary differential equation}
We proceed to solve the Cauchy problem 
$(\ref{eq:third})$.
To this end,
we have to evaluate the solutions of the $5^{\rm th}$--degree characteristic equation in the variable $Z$
\begin{equation}\label{eq:fchar}
Z^{5}+\eta Z^{4}+2\la_{n}Z^{3}+\la_{n} (2\eta-\be)Z^{2}+(\la_{n}^2-ab)Z+\la_{n}^2(\eta-\be)-\eta ab=0\,.
\end{equation}
By means of the singular perturbation theory we get  the five solutions of \eqref{eq:fchar}: one is a real number $r_n$  and the other four $i\omega_n$, $-i\overline{\omega_n}$, $i\zeta_n$, $-i\overline{\zeta_n}$ are pairwise complex conjugate numbers. Moreover, $r_n$, $\omega_n$ and $\zeta_n$ exhibit the following
asymptotic behavior as $n$ tends to $\infty$: 
\begin{equation}\label{eq:lambda1}
r_{n}=\be-\eta
-{\beta\big(\be-\eta\big)^2\over\la_{n}}+O\Big({1\over{\la_{n}^{2}}}\Big)
=\be-\eta
+O\Big({1\over{\la_{n}}}\Big)
\,,
\end{equation}
\begin{multline}\label{eq:lambda2}
\omega_{n}=
\sqrt{\la_{n}}+{\be\over2}\Big({3\over4}\beta-\eta\Big){1\over\sqrt{\lambda_{n}}}
+i \Big[{\be\over 2}
-\Big({\beta\big(\be-\eta\big)^2\over2}+\frac{ab}{2\beta}\Big){1\over\la_{n}} \Big]
+O\Big({1\over{\la_{n}^{3/2}}}\Big)
\\
=
\sqrt{\la_{n}}+i{\be\over 2}
+O\Big({1\over{\sqrt{\la_{n}}}}\Big)
\,,
\end{multline}
\begin{equation}\label{eq:lambda4}
\zeta_{n}
=
\sqrt{\la_{n}}+{\eta ab\over 2\beta\la_{n}^{3/2}}
+i\Big({ab\over 2\beta\la_{n}}+{a^2b^2\over 2\beta^3\la^2_{n}}\Big)
+O\Big({1\over{\la_{n}^{5/2}}}\Big)
=\sqrt{\la_{n}}+i{ab\over 2\beta\la_{n}}+O\Big({1\over{\la_{n}^{3/2}}}\Big)
\,.
\end{equation}
Therefore, we are able to write the solution  $f_{1n}(t)$ of (\ref{eq:third}) in the form
\begin{equation}\label{eq:f1j}
f_{1n}(t)=R_{n}e^{r_{n} t}+C_{n}e^{i\omega_{n} t}+\overline{C_{n}}e^{-i\overline{\omega_{n}}t}
+D_{n}e^{i\zeta_{n} t}+\overline{D_{n}}e^{-i\overline{\zeta_{n}}t}
\,,
\end{equation}
where the coefficients $R_{n}\in\R$ and $C_{n},D_{n}\in\C$ are unknown.
Since the function $f_{1n}(t)$ have to satisfy the initial conditions in $(\ref{eq:third})$, 
to determine $R_{n}$, $C_{n}$ and $D_{n}$ 
we will  solve 
 the system
\begin{equation}\label{vandermonde}
\left \{\begin{array}{l}
R_{n}+C_{n}+\overline{C_{n}}+D_{n}+\overline{D_{n}}=f_{1n}(0),\\ 
\\
r_{n} R_{n}+i\omega_{n} C_{n}-i\overline{\omega_{n}C_{n}}+i\zeta_{n} D_{n}-i\overline{\zeta_{n}D_{n}}
=f_{1n}'(0),\\
\\
r_{n}^2R_{n}-\omega_{n}^2C_{n}-\overline{\omega_{n}^2C_{n}}-\zeta_{n}^2D_{n}-\overline{\zeta_{n}^2D_{n}}=f_{1n}''(0),\\
\\
r_{n}^3R_{n}-i \omega_{n}^3C_{n}+i\overline{\omega_{n} ^3C_{n}}-i\zeta_{n}^3D_{n}+i\overline{\zeta_{n}^3D_{n}}=f_{1n}'''(0),\\
\\
r_{n}^4R_{n}+ \omega_{n} ^4C_{n}+\overline{\omega_{n}^4C_{n}}+\zeta_{n}^4D_{n}+\overline{\zeta_{n}^4D_{n}}=f_{1n}^{(4)}(0).
\end{array}\right .
\end{equation}
Indeed, we obtain that the coefficients 
have the following asymptotic behavior as $n$ tends to $\infty$:
\begin{equation}\label{eq:asy_R}
R_n={\beta\over\la_n}(\alpha_{1n}(\beta-\eta)+\rho_{1n})
+( \alpha_{1n}+\rho_{1n}+\alpha_{2n}+\rho_{2n})O\Big({1\over{\la_n^{2}}}\Big),
\end{equation}
\begin{multline}\label{eq:asy_C}
C_{n}={\alpha_{1n}\over 2}
-\frac{i}{4\beta}
\big(\beta^2\alpha_{1n}+2\beta\rho_{1n}
+2a\alpha_{2n}\Big)
\frac{1}{\la_{n}^{1/2}}
+\frac{1}{2\beta^2}
\big((ab-\beta^3(\beta-\eta))\alpha_{1n}
\\
-\beta(\beta^2\rho_{1n}+\eta a\alpha_{2n})
-\beta a\rho_{2n}\Big)\frac1{\la_{n}}
+(\alpha_{1n}+\rho_{1n}+\alpha_{2n}+\rho_{2n})O\Big({1\over{\la_{n}^{3/2}}}\Big)
\end{multline}
\begin{multline}\label{eq:asy_D}
D_n=i\frac{a\alpha_{2n}}{2\beta\la_{n}^{1/2}}
+\frac{a}{2\beta^2}
\big(\beta\eta\alpha_{2n}+\beta\rho_{2n}
-b\alpha_{1n}\big)
\frac1{\la_n}
+\frac{i}{2\beta^3}
\big(2a^2b\alpha_{2n}-\eta\beta^2a\rho_{2n}+2\eta\beta ab\alpha_{1n}+\beta ab\rho_{1n}\big)
\frac1{\la_n^{3/2}}
\\
+(\alpha_{1n}+\rho_{1n}+ \alpha_{2n}+\rho_{2n})O\Big({1\over{\la_n^{2}}}\Big)
.
\end{multline}
Accordingly, we can write $f_{1n}(t)$
by means of formula \eqref{eq:f1j}, where the coefficients $R_n$, $C_n$ and $D_n$ are given by formulas \eqref{eq:asy_R}-\eqref{eq:asy_D} respectively. Moreover, 
thanks to \eqref{eq:f2j}, we can also get the expression 
for $f_{2n}(t)$, that is
\begin{equation}\label{eq:f2n}
f_{2n}(t)
=\Upsilon_n\big(R_{n}e^{r_{n} t}+C_{n}e^{i\omega_{n} t}+\overline{C_{n}}e^{-i\overline{\omega_{n}}t}+D_{n}e^{i\zeta_{n} t}+\overline{D_{n}}e^{-i\overline{\zeta_{n}}t}\big)
\,.
\end{equation}
We will observe that the function $f_{2n}(t)$ can be written in a more handleable form.
To this end, first we recall the following result (see e.g. \cite[Section 6]{LoretiSforza1})
\begin{lemma}\label{eq:appr-ce}
Approximated solutions of the cubic equation
\begin{equation}\label{eq:char}
Z^{3}+\eta Z^{2}+\la_{n}Z+\la_{n}(\eta-\be)=0\,,
\end{equation}
are given by
\begin{equation}\label{eq:lambda10}
r_{n}=\be-\eta
-{\beta\big(\be-\eta\big)^2\over\la_{n}}+O\Big({1\over{\la_{n}^{2}}}\Big)
\,,
\end{equation}
\begin{equation}\label{eq:lambda20}
z_{n}=
-{\be\over 2}
+{\beta\big(\be-\eta\big)^2\over2}{1\over\la_{n}} 
+i\Big[\sqrt{\la_{n}}+{\be\over2}\Big({3\over4}\beta-\eta\Big){1\over\sqrt{\lambda_{n}}}\Big]
+O\Big({1\over{\la_{n}^{3/2}}}\Big)
\,.
\end{equation}
\end{lemma}
Therefore, comparing \eqref{eq:lambda1} with \eqref{eq:lambda10}, we have that the numbers $r_n$
are approximated solutions of \eqref{eq:char}, and hence
the function $t\to R_{n}e^{r_{n} t}$
is a solution of the third order differential equation
\begin{equation}\label{eq:=0}
g{'''}+\eta g{''}+\la_n g{'}+\la_n (\eta-\be)g=0
\qquad \mbox{in}\,\,\,\R\,.
\end{equation}
\begin{lemma}\label{eq:app-sol}
The numbers $i\omega_n$, with $\omega_n$ defined by \eqref{eq:lambda2}, 
are approximated solutions of the cubic equation
\begin{equation*}
Z^{3}+\eta Z^{2}+\la_{n}Z+\la_{n}(\eta-\be)=-\frac{ab}{\beta}\,.
\end{equation*}
\end{lemma}
\begin{Proof}
The comparison of \eqref{eq:lambda2} with \eqref{eq:lambda20} yields
\begin{equation*}
i\omega_n=z_n+\frac{ab}{2\beta\lambda_n}
\,.
\end{equation*}
Since
\begin{multline*}
(i\omega_n)^{3}+\eta (i\omega_n)^{2}+\la_{n}i\omega_n+\la_{n}(\eta-\be)
\\
=z_n^3+\eta z_n^2+\la_{n}z_n+\la_{n}(\eta-\be)
+3z_n^2\frac{ab}{2\beta\lambda_n}+3z_n\frac{a^2b^2}{4\beta^2\lambda_n^2}+\frac{a^3b^3}{8\beta^3\lambda_n^3}
+2\eta z_n\frac{ab}{2\beta\lambda_n}+\eta\frac{a^2b^2}{4\beta^2\lambda_n^2}
+\frac{ab}{2\beta}
\,,
\end{multline*}
and in virtue  of Lemma \ref{eq:appr-ce} we have
\begin{equation*}
z_n^3+\eta z_n^2+\la_{n}z_n+\la_{n}(\eta-\be)=0,
\end{equation*}
then we get
\begin{equation*}
(i\omega_n)^{3}+\eta (i\omega_n)^{2}+\la_{n}i\omega_n+\la_{n}(\eta-\be)
=-\frac{3ab}{2\beta}+\frac{ab}{2\beta}+O\Big({1\over{\sqrt{\la_{n}}}}\Big)=-\frac{ab}{\beta}+O\Big({1\over{\sqrt{\la_{n}}}}\Big)
\,.
\end{equation*}
that is, our claim holds true.
\end{Proof}
Thanks to Lemma \ref{eq:app-sol}, the numbers $i\omega_n$ and  their conjugate numbers $-i\overline{\omega_n}$ are approximated solutions of the cubic equation
\begin{equation*}
Z^{3}+\eta Z^{2}+\la_{n}Z+\la_{n}(\eta-\be)=-\frac{ab}{\beta}\,,
\end{equation*}
so, it follows that the function
$
t\to C_{n}e^{i\omega_{n} t}+\overline{C_{n}}e^{-i\overline{\omega_{n}}t}
$
is a solution of the third order differential equation
\begin{equation}\label{eq:not=0}
g{'''}+\eta g{''}+\la_n g{'}+\la_n (\eta-\be)g=-\frac{ab}{\beta}g
\qquad \mbox{in}\,\,\,\R\,.
\end{equation}
In virtue of \eqref{eq:=0} and \eqref{eq:not=0}, the function
$$
g_n(t)= R_{n}e^{r_{n} t}+C_{n}e^{i\omega_{n} t}+\overline{C_{n}}e^{-i\overline{\omega_{n}}t}
$$
is a solution of the third order differential equation
\begin{equation}
g{'''}+\eta g{''}+\la_n g{'}+\la_n (\eta-\be)g=-\frac{ab}{\beta}(C_{n}e^{i\omega_{n} t}+\overline{C_{n}}e^{-i\overline{\omega_{n}}t})
\qquad \mbox{in}\,\,\,\R\,.
\end{equation}
Therefore, we can apply Lemma \ref{le:third} with $h(t)=-\frac{ab}{\beta}(C_{n}e^{i\omega_{n} t}+\overline{C_{n}}e^{-i\overline{\omega_{n}}t})$:
thanks to \eqref{eq:second} and \eqref{eq:f2j'}, we have
\begin{equation}\label{eq:C_{n}e}
\Upsilon_n(g_n(t))
=-\frac1ae^{-\eta t}\big(g_n''(0)+\la_n g_n(0)\big)
+\frac{b}{\beta} \int_0^t e^{-\eta(t-s)}(C_{n}e^{i\omega_{n} s}+\overline{C_{n}}e^{-i\overline{\omega_{n}}s})ds\,.
\end{equation}
From \eqref{vandermonde} and \eqref{eq:third} it follows that
\begin{equation*}
\begin{split}
  g_n''(0)
  &=f''_{1n}(0)+\zeta_{n}^2D_{n}+\overline{\zeta_{n}^2D_{n}}
  =-\la_{n}\alpha_{1n}-a\alpha_{2n}+\zeta_{n}^2D_{n}+\overline{\zeta_{n}^2D_{n}}
  \\
  \lambda_n g_n(0)
  &= \lambda_nf_{1n}(0)-\lambda_n D_{n}-\lambda_n\overline{D_{n}}
  =\la_{n}\alpha_{1n}-\lambda_n D_{n}-\lambda_n\overline{D_{n}}
  \,.
\end{split}
\end{equation*}
Thanks to \eqref{eq:lambda4} we have
$\zeta_{n}^2-\lambda_n=O\Big({1\over{\sqrt{\la_{n}}}}\Big)$,
so we see that
\begin{equation*}
g_n''(0)+\la_n g_n(0)=
-a\alpha_{2n}+( \alpha_{1n}+\rho_{1n}+\alpha_{2n}+\rho_{2n})O\Big({1\over{\la_n}}\Big)
\,.
\end{equation*}
Moreover
\begin{equation*}
\int_0^t e^{-\eta(t-s)}e^{i\omega_{n} s}ds
=\frac{1}{\eta+i\omega_n}\big(e^{i\omega_{n} t}-e^{-\eta t}\big)\,.
\end{equation*}
Set
\begin{equation}
c_n=\frac{b}{\beta(\eta+i\omega_n)}
\,,
\end{equation}
from  \eqref{eq:C_{n}e} we obtain
\begin{equation}\label{eq:R_{n}+C_{n}e}
\Upsilon_n(R_{n}e^{r_{n} t}+C_{n}e^{i\omega_{n} t}+\overline{C_{n}}e^{-i\overline{\omega_{n}}t})
=
c_nC_{n}e^{i\omega_{n} t}
+\overline{c_nC_{n}}e^{-i\overline{\omega_{n}}t}
+\big(\alpha_{2n}
-2\Re(c_nC_{n})\big)e^{-\eta t}\,.
\end{equation}
Moreover, thanks to \eqref{eq:Yexp} we have
\begin{equation*}
\Upsilon_n(e^{i\zeta_{n} t})
=
\frac1a
\Big(\zeta_{n}^2-\la_n+\frac{\beta\la_n}{\eta+i\zeta_{n}}\Big)e^{i\zeta_{n} t}-\frac{\beta\la_n}{a(\eta+i\zeta_{n})}e^{-\eta t}
\,.
\end{equation*}
Therefore, 
if we define
\begin{equation}
d_n=\frac1a
\Big(\zeta_{n}^2-\la_n+\frac{\beta\la_n}{\eta+i\zeta_{n}}\Big)
\,,
\end{equation}
and
\begin{equation}
E_n=\alpha_{2n}
-2\Re(c_nC_{n})-\frac{2\beta\la_n}{a}\Re\bigg(\frac{D_{n}}{\eta+i\zeta_{n}}\bigg)\,,
\end{equation}
thanks to \eqref{eq:f2n} and \eqref{eq:R_{n}+C_{n}e},  $f_{2n}(t)$ can be written in the following form
\begin{equation}\label{eq:f2jbis0}
f_{2n}(t)
=
d_nD_{n}e^{i\zeta_{n} t}+\overline{d_nD_{n}}e^{-i\overline{\zeta_{n}}t}
+c_nC_{n}e^{i\omega_{n} t}+\overline{c_nC_{n}}e^{-i\overline{\omega_{n}}t}
+E_ne^{-\eta t}
\,.
\end{equation}
We also note that
\begin{equation}\label{eq:cndn0}
|d_n|\asymp|\zeta_n|\asymp \sqrt{\la_{n}}\,,
\qquad
|c_n|\le\frac{M}{|\omega_n|}\,.
\end{equation}
The proof of the following lemma is straightforward in virtue of 
\eqref{eq:asy_D} and \eqref{eq:cndn0}, so we omit it.
\begin{lemma}\label{le:mathcalE0}
Set
\begin{equation*}
E_n=\alpha_{2n}
-2\Re(c_nC_{n})-\frac{2\beta\la_n}{a}\Re\bigg(\frac{D_{n}}{\eta+i\zeta_{n}}\bigg)\,,
\end{equation*}
there exists a constant $M>0$ such that
\begin{equation*}
\Big|\sum_{ n= 1}^{\infty} E_n\Big|^2\le 
M \sum_{ n= 1}^{\infty}
\Big(|C_{n}|^2+ |d_nD_{n}|^2\Big)
\,.
\end{equation*}
\end{lemma}

Now, we state and prove some properties about the coefficients, that show some differences with respect to  the analogous ones in \cite{LoretiSforza1, LoretiSforza3}.

\begin{lemma}
The following statements hold true.
\begin{itemize}
\item[(i)] For any $n\in\N$ one has
\begin{equation}\label{eq:|Cj2|}
|C_{n}|^2+\lambda_n |D_{n}|^2
\asymp
\frac{1}{\la_{n}}
\big(\alpha^2_{1n}\la_{n}+\rho^2_{1n}
+\alpha^2_{2n}\la_{n}+\rho^2_{2n}\big).
\end{equation}
\item[(ii)] There exists a constant $M>0$ such that  for any $n\in\N$ one has
\begin{equation}\label{eq:C1overC2}
\vert R_{n}\vert\le {M \over {{\lambda^{1/2}_n}}}\Big(|C_{n}|^2+\lambda_n |D_{n}|^2\Big)^{1/2}\,.  
\end{equation}

\end{itemize}
\end{lemma}
\begin{Proof}
(i) From \eqref{eq:asy_C} it follows that
\begin{multline}\label{eq:modCn}
|C_{n}|^2
={1\over 4}\alpha_{1n}^2
+\frac{1}{16\beta^2}
\big(\beta^2\alpha_{1n}+2\beta\rho_{1n}
+2a\alpha_{2n}\Big)^2
\frac{1}{\la_{n}}
\\
+\frac{\alpha_{1n}}{2\beta^2}
\big((ab-\beta^3(\beta-\eta))\alpha_{1n}
-\beta(\beta^2\rho_{1n}+\eta a\alpha_{2n})
-\beta a\rho_{2n}\Big)\frac1{\la_{n}}
\\
+(\alpha_{1n}^2+\rho_{1n}^2+\alpha_{2n}^2+\rho_{2n}^2)O\Big({1\over{\la_{n}^{2}}}\Big)
\,.
\end{multline}
Moreover, from \eqref{eq:asy_D} we deduce that
\begin{multline*}
\la_{n}^{1/2}D_n=i\frac{a\alpha_{2n}}{2\beta}
+\frac{a}{2\beta^2}
\big(\beta\eta\alpha_{2n}+\beta\rho_{2n}
-b\alpha_{1n}\big)
\frac1{\la_{n}^{1/2}}
\\
+\frac{i}{2\beta^3}
\big(2a^2b\alpha_{2n}+2\eta\beta ab\alpha_{1n}+\beta ab\rho_{1n}-\eta\beta^2a\rho_{2n}\big)
\frac1{\la_n}
+(\alpha_{1n}+\rho_{1n}+ \alpha_{2n}+\rho_{2n})O\Big({1\over{\la_n^{3/2}}}\Big)
,
\end{multline*}
whence
\begin{multline}\label{eq:modDn}
\la_{n}|D_n|^2=\frac{a^2\alpha_{2n}^2}{4\beta^2}
+\frac{a^2}{4\beta^4}
\big(\beta\eta\alpha_{2n}+\beta\rho_{2n}
-b\alpha_{1n}\big)^2
\frac1{\la_{n}}
\\
+\frac{a\alpha_{2n}}{2\beta^4}
\big(2a^2b\alpha_{2n}+2\eta\beta ab\alpha_{1n}+\beta ab\rho_{1n}-\eta\beta^2a\rho_{2n}\big)
\frac1{\la_n}
+(\alpha_{1n}^2+\rho_{1n}^2+\alpha_{2n}^2+\rho_{2n}^2)O\Big({1\over{\la_n^{2}}}\Big)
.
\end{multline}
Now, putting together \eqref{eq:modCn} and \eqref{eq:modDn}, we have
\begin{multline*}
|C_{n}|^2+\la_{n}|D_n|^2
=\frac14
\Big(\alpha^2_{1n}+\frac{\rho^2_{1n}}{\la_{n}}
+\frac{a^2}{\beta^2}\Big(
\alpha^2_{2n}
+\frac{\rho^2_{2n}}{\la_{n}}\Big)\Big)
\\
+\frac{1}{16\beta^2}
\big(\beta^2\alpha_{1n}+2a\alpha_{2n}\big)^2
\frac{1}{\la_{n}}
+\frac{\rho_{1n}}{4\beta}\big(\beta^2\alpha_{1n}+2a\alpha_{2n}\big)
\frac{1}{\la_{n}}
\\
+\frac{\alpha_{1n}}{2\beta^2}
\big((ab-\beta^3(\beta-\eta))\alpha_{1n}
-\beta(\beta^2\rho_{1n}+\eta a\alpha_{2n})
-\beta a\rho_{2n}\Big)\frac1{\la_{n}}
\\
+\frac{a^2}{4\beta^4}
\big(\beta\eta\alpha_{2n}
-b\alpha_{1n}\big)^2
\frac1{\la_{n}}
+\frac{a^2\rho_{2n}}{2\beta^3}
\big(\beta\eta\alpha_{2n}
-b\alpha_{1n}\big)
\frac1{\la_{n}}
\\
+\frac{a\alpha_{2n}}{2\beta^4}
\big(2a^2b\alpha_{2n}+2\eta\beta ab\alpha_{1n}+\beta ab\rho_{1n}-\eta\beta^2a\rho_{2n}\big)
\frac1{\la_n}
+\big(\alpha^2_{1n}+\rho^2_{1n}
+\alpha^2_{2n}+\rho^2_{2n}
\big)O\Big({1\over{\la_{n}^{2}}}\Big)
\,.
\end{multline*}
We can neglect the indices $n\in\N$ such that $\alpha_{1n}=\rho_{1n}=\alpha_{2n}=\rho_{2n}=0$, because the present evaluation will be used in summing series. So, we can assume that for any $n\in\N$ $(\alpha_{1n},\rho_{1n},\alpha_{2n},\rho_{2n})\not=(0,0,0,0)$, and hence by the previous formula
we obtain
\begin{multline*}
\frac{|C_{n}|^2+\la_{n}|D_n|^2}{\alpha^2_{1n}+\frac{\rho^2_{1n}}{\la_{n}}
+\frac{a^2}{\beta^2}\Big(
\alpha^2_{2n}
+\frac{\rho^2_{2n}}{\la_{n}}\Big)}
\\
={1\over 4}
+\frac{\big(
\alpha_{1n}^2+(\alpha_{1n}+\alpha_{2n})(\rho_{1n}+\alpha_{2n}+\rho_{2n})
\big)O\Big({1\over{\la_{n}}}\Big)
}
{\alpha^2_{1n}+\frac{\rho^2_{1n}}{\la_{n}}
+\frac{a^2}{\beta^2}\Big(\alpha^2_{2n}
+\frac{\rho^2_{2n}}{\la_{n}}\Big)}
\to\frac14\,,\qquad\mbox{as}\qquad n\to\infty\,,
\end{multline*}
taking into account, for example, that 
$$
{\a_{1n}\rho_{1n}\over{\la_{n}}}
={\a_{1n}\over{\la^{1/3}_{n}}}
{\rho_{1n}\over{\la^{2/3}_{n}}}
\le{\a^2_{1n}\over{\la^{2/3}_{n}}}+{\rho_{1n}^2\over{\la^{4/3}_{n}}}
\,.
$$
In conclusion, \eqref{eq:|Cj2|} holds true.

\noindent (ii)
From \eqref{eq:asy_R} we have
\begin{equation*}
|R_n|^2={\beta^2\over\la_n^2}\big(\alpha_{1n}(\beta-\eta)+\rho_{1n}\big)^2
+\big( \alpha_{1n}+\rho_{1n}\big)\big( \alpha_{1n}+\rho_{1n}+\alpha_{2n}+\rho_{2n}\big)O\Big({1\over{\la_n^3}}\Big).
\end{equation*}
Moreover, thanks to \eqref{eq:|Cj2|}, there exists a constant $c>0$ such that
\begin{equation*}
|C_{n}|^2+\lambda_n |D_{n}|^2\ge \frac{c}{\la_n}\big(\alpha^2_{1n}\la_{n}+\rho^2_{1n}
+\alpha^2_{2n}\la_{n}+\rho^2_{2n}\big)\,.
\end{equation*}
Therefore, from the above formulas we get
\begin{equation*}
\frac{|R_n|^2}{|C_{n}|^2+\lambda_n |D_{n}|^2}
\le{1\over c\la_n}
\frac{\beta^2(\alpha_{1n}(\beta-\eta)+\rho_{1n})^2
+( \alpha_{1n}+\rho_{1n})( \alpha_{1n}+\rho_{1n}+\alpha_{2n}+\rho_{2n})O\Big({1\over{\la_n}}\Big)}
{\alpha^2_{1n}\la_{n}+\rho^2_{1n}+\alpha^2_{2n}\la_{n}+\rho^2_{2n}}
\,,
\end{equation*}
that is, \eqref{eq:C1overC2} follows.
\end{Proof}
In conclusion, taking into account of any result of the present section we have  proved the following representation formula for the solution of the coupled system.
\begin{theorem}\label{th:repres}
The solution of problem \eqref{eq:system} can be written as series in the following way
\begin{equation}\label{eq:u1}
\begin{split}
u_1(t) &=\sum_{n=1}^{\infty}\Big(C_{n}e^{i\omega_{n} t}+\overline{C_{n}}e^{-i\overline{\omega_{n}}t}
+R_{n}e^{r_{n} t}+D_{n}e^{i\zeta_{n} t}+\overline{D_{n}}e^{-i\overline{\zeta_{n}}t}\Big)w_{n}
\,,
\\
u_2(t) &=\sum_{n=1}^{\infty}
\Big(d_nD_{n}e^{i\zeta_{n} t}+\overline{d_nD_{n}}e^{-i\overline{\zeta_{n}}t}
+c_nC_{n}e^{i\omega_{n} t}+\overline{c_nC_{n}}e^{-i\overline{\omega_{n}}t}+E_n e^{-\eta t}
\Big)w_{n}
\,,
\end{split}
\end{equation}
where
\begin{equation*}
r_{n}
=\be-\eta
+O\Big({1\over{\la_{n}}}\Big)
\,,
\end{equation*}
\begin{equation*}
\omega_{n}
=
\sqrt{\la_{n}}+i{\be\over 2}
+O\Big({1\over{\sqrt{\la_{n}}}}\Big)
\,,
\end{equation*}
\begin{equation*}
\zeta_{n}
=\sqrt{\la_{n}}+i{ab\over 2\beta\la_{n}}+O\Big({1\over{\la_{n}^{3/2}}}\Big)
\,,
\end{equation*}
\begin{equation*}
\vert R_{n}\vert\le {M \over {{\lambda^{1/2}_n}}}\Big(|C_{n}|^2+|d_nD_{n}|^2\Big)^{1/2}\,,
\quad 
\Big|\sum_{ n= 1}^{\infty} E_n\Big|^2\le 
M \sum_{ n= 1}^{\infty}
\Big(|C_{n}|^2+ |d_nD_{n}|^2\Big)
\,,
\end{equation*}
\begin{equation*}
|d_n|\asymp \sqrt{\la_{n}}\,,
\quad
|c_n|\le\frac{M}{\sqrt{\la_n}}\,,
\qquad (M>0)
\end{equation*}
\begin{equation*}
\sum_{ n= 1}^{\infty}\la_{n}\Big(|C_{n}|^2+|d_nD_{n}|^2\Big)
\asymp
\|u_1^0\|^2_{D(\sqrt{L})}+\|u_1^1\|^2_{H}+\|u_2^0\|^2_{D(\sqrt{L})}+\|u_2^1\|^2_{H}
\,.
\end{equation*}
\end{theorem}
\section{Ingham type estimates}\label{se:invdir}
Our goal is to prove an inverse inequality and a direct inequality for the pair $(u_1,u_2)$ defined by
\begin{equation}\label{eq:vsum1}
\begin{split}
u_1(t) &=\sum_{n=1}^{\infty}\Big(C_{n}e^{i\omega_{n} t}+\overline{C_{n}}e^{-i\overline{\omega_{n}}t}
+R_{n}e^{r_{n} t}+D_{n}e^{i\zeta_{n} t}+\overline{D_{n}}e^{-i\overline{\zeta_{n}}t}\Big)
\,,
\\
u_2(t) &=\sum_{n=1}^{\infty}
\Big(d_nD_{n}e^{i\zeta_{n} t}+\overline{d_nD_{n}}e^{-i\overline{\zeta_{n}}t}
+c_nC_{n}e^{i\omega_{n} t}+\overline{c_nC_{n}}e^{-i\overline{\omega_{n}}t}\Big)
+\mathcal {E} e^{-\eta t}
\,,
\end{split}
\end{equation}
with $\om_{n}\,,C_{n}\,,\zeta_{n}\,,D_{n}, d_n,c_n\in\C$
and
$r_{n}\,,R_{n}\,,\mathcal {E}\in\R$. 
We will assume that  there exist
$\gamma>0$, $\alpha,\chi\in\R$, $n'\in\N$, $\mu>0$, $\nu> 1/2$, 
such that
\begin{equation}\label{eq:hom1}
\liminf_{n\to\infty}({\Re}\om_{n+1}-{\Re}\om_{n})=\liminf_{n\to\infty}({\Re}\zeta_{n+1}-{\Re} \zeta_{n})=\gamma\,,
\end{equation}
\begin{equation}\label{eq:hom2}
\begin{split}
\lim_{n\to\infty}{\Im}\om_n&=\alpha>0
\,,
\\
\lim_{n\to\infty}r_n&=\chi<0\,,
\\
\lim_{n\to\infty}\Im \zeta_{n}&=0\,,
\end{split}
\end{equation}
\begin{equation}\label{eq:cndn}
|d_n|\asymp|\zeta_n|
\,,
\qquad
|c_n|\le\frac{M}{|\omega_n|}\,,
\end{equation}
\begin{equation}\label{eq:hom3}
|R_n|\le \frac{\mu}{n^{\nu}}\Big(|C_{n}|^2+ |d_nD_{n}|^2\Big)^{1/2}\,\quad\forall\ n\ge n'\,,
\qquad
|R_n|\le \mu\Big(|C_{n}|^2+ |d_nD_{n}|^2\Big)^{1/2}\,\quad\forall\ n\le n'\,.
\end{equation}
\subsection{Outline of the proof}\label{se:outline}
Before to proceed with our computations, we will outline briefly our reasoning. Firstly, to shorten our formulas we introduce the following notations
\begin{equation}\label{eq:f11,f12}
{\mathcal U}_{1}^C(t)=\sum_{n=1}^{\infty}\big(C_{n}e^{i\omega_{n} t}+\overline{C_{n}}e^{-i\overline{\omega_{n}}t}\big),
\quad
{\mathcal U}_{1}^D(t)=\sum_{n=1}^{\infty}\big(D_{n}e^{i\zeta_{n} t}+\overline{D_{n}}e^{-i\overline{\zeta_{n}}t}\big),
\quad
{\mathcal U}_{1}^R(t)=\sum_{n=1}^{\infty} R_{n}e^{r_{n} t},
\end{equation}
\begin{equation}\label{eq:f22,f21}
{\mathcal U}_{2}^D(t)=\sum_{n=1}^\infty \big(d_nD_{n}e^{i \zeta_{n}t}+\overline{d_nD_{n}}e^{-i \overline{\zeta_n}t}\big),
\qquad
{\mathcal U}_{2}^C(t)=\sum_{n=1}^\infty \big(c_nC_{n}e^{i\om_nt}+\overline{c_nC_{n}}e^{-i\overline{\omega_n}t}\big)\,,
\end{equation}
so  we can write the functions $u_1$, $u_2$ as 
\begin{equation*}
u_1={\mathcal U}_{1}^C+{\mathcal U}_{1}^D+{\mathcal U}_{1}^R,
\qquad
u_2-\mathcal {E} e^{-\eta t}={\mathcal U}_{2}^D+{\mathcal U}_{2}^C.
\end{equation*}
If $k(t)$  is a suitable positive function, see (\ref{eq:k}) below, our first goal will be to estimate 
\begin{equation*}
\int_{0}^{\infty} k(t)| {\mathcal U}_{1}^C(t)+{\mathcal U}_{1}^D(t)+{\mathcal U}_{1}^R(t)|^2\ dt
+\int_{0}^{\infty} k(t)| {\mathcal U}_{2}^D(t)+ {\mathcal U}_{2}^C(t)|^2\ dt
\,,
\end{equation*}
unless a finite number of terms in the series.

By reason of $2 ab\ge-\frac12a^2-2b^2$ we have $|a+b|^2\ge\frac12a^2-b^2$, so we can observe that
\begin{equation*}
|{\mathcal U}_{1}^C(t)+{\mathcal U}_{1}^D(t)+{\mathcal U}_{1}^R(t)|^2
\ge
\frac12|{\mathcal U}_{1}^C(t)|^2-|{\mathcal U}_{1}^D(t)+{\mathcal U}_{1}^R(t)|^2
\ge
\frac12|{\mathcal U}_{1}^C(t)|^2-2|{\mathcal U}_{1}^D(t)|^2-2|{\mathcal U}_{1}^R(t)|^2
\,,
\end{equation*}
\begin{equation*}
| {\mathcal U}_{2}^D(t)+ {\mathcal U}_{2}^C(t)|^2
\ge
\frac12| {\mathcal U}_{2}^D(t)|^2-|{\mathcal U}_{2}^C(t)|^2\,.
\end{equation*}
Bearing in mind \eqref{eq:hom3}, since $k(t)$ is positive from the above  inequalities we can deduce
\begin{multline*}
\int_{0}^{\infty} k(t)| {\mathcal U}_{1}^C(t)+{\mathcal U}_{1}^D(t)+{\mathcal U}_{1}^R(t)|^2\ dt
+\int_{0}^{\infty} k(t)| {\mathcal U}_{2}^D(t)+ {\mathcal U}_{2}^C(t)|^2\ dt
\\
\ge
\int_{0}^{\infty} k(t) \Big(\frac12| {\mathcal U}_{1}^C(t)|^2-2|{\mathcal U}_{1}^D(t)|^2\Big)\ d t
+\int_{0}^{\infty} k(t) \Big(\frac12| {\mathcal U}_{2}^D(t)|^2-|{\mathcal U}_{2}^C(t)|^2\Big)\ d t
\\
-2\int_{0}^{\infty} k(t) | {\mathcal U}_{1}^R(t)|^2\ d t
\,.
\end{multline*}
In virtue of \eqref{eq:cndn} we can control the term $\int_{0}^{\infty} k(t) {\mathcal U}_{1}^D(t) d t$ (resp. $\int_{0}^{\infty} k(t) {\mathcal U}_{2}^C(t) d t$) by means of 
\break
$\int_{0}^{\infty} k(t) {\mathcal U}_{2}^D(t) d t$ 
(resp. $\int_{0}^{\infty} k(t) {\mathcal U}_{1}^C(t)\d t$). Therefore, it is convenient to write the previous formula in the following way
\begin{multline}\label{eq:f1+f2}
\int_{0}^{\infty} k(t)| {\mathcal U}_{1}^C(t)+{\mathcal U}_{1}^D(t)+{\mathcal U}_{1}^R(t)|^2\ dt
+\int_{0}^{\infty} k(t)| {\mathcal U}_{2}^D(t)+ {\mathcal U}_{2}^C(t)|^2\ dt
\\
\ge
\frac12\int_{0}^{\infty} k(t) \Big(| {\mathcal U}_{1}^C(t)|^2-2|{\mathcal U}_{2}^C(t)|^2\Big)\ d t
+\frac12\int_{0}^{\infty} k(t) \Big(| {\mathcal U}_{2}^D(t)|^2-4|{\mathcal U}_{1}^D(t)|^2\Big)\ d t
\\
-2\int_{0}^{\infty} k(t) | {\mathcal U}_{1}^R(t)|^2\ d t
\,.
\end{multline}
We will give a lower bound estimate for  $\int_{0}^{\infty} k(t) | {\mathcal U}_{1}^C(t)|^2 d t $ and 
$\int_{0}^{\infty} k(t) | {\mathcal U}_{2}^D(t)|^2 d t $, and, on the contrary,  an upper bound estimate for $\int_{0}^{\infty} k(t) | {\mathcal U}_{2}^C(t)|^2 d t $,
$\int_{0}^{\infty} k(t) | {\mathcal U}_{1}^D(t)|^2 d t $ and 
$\int_{0}^{\infty} k(t) | {\mathcal U}_{1}^R(t)|^2 d t $.  So, thanks to \eqref{eq:f1+f2}, we will be able to prove an inverse estimate.

Moreover, if we will assume an additional condition on the coefficients of the series, we will be able to prove an inverse inequality with a better estimate for the control time.
Indeed, the additional assumption will allow us to control all terms 
$\int_{0}^{\infty} k(t) |{\mathcal U}_{1}^D(t)|^2 d t$, 
$\int_{0}^{\infty} k(t) |{\mathcal U}_{2}^C(t)|^2d t$ and
$\int_{0}^{\infty} k(t) | {\mathcal U}_{1}^R(t)|^2d t$ 
by means of 
$\int_{0}^{\infty} k(t) | {\mathcal U}_{2}^D(t)|^2 d t$. 
In this way the estimate of the term $\int_{0}^{\infty} k(t) | {\mathcal U}_{1}^C(t)|^2\ d t$ can be done with the help of an idea used previously in \cite{LoretiSforza1}.
In fact in this case we will use the following inequality
\begin{multline}\label{eq:f1+f2bis}
\int_{0}^{\infty} k(t)| {\mathcal U}_{1}^C(t)+{\mathcal U}_{1}^D(t)+{\mathcal U}_{1}^R(t)|^2\ dt
+\int_{0}^{\infty} k(t)| {\mathcal U}_{2}^D(t)+ {\mathcal U}_{2}^C(t)|^2\ dt
\\
\ge
\frac12\int_{0}^{\infty} k(t) | {\mathcal U}_{1}^C(t)|^2\ d t
+\frac12\int_{0}^{\infty} k(t) \Big(| {\mathcal U}_{2}^D(t)|^2-4|{\mathcal U}_{1}^D(t)|^2-2|{\mathcal U}_{2}^C(t)|^2-4| {\mathcal U}_{1}^R(t)|^2
\Big)\ d t
\,.
\end{multline}

\subsection{Technical results}
In order to avoid repetitions and simplify the proofs of the main theorems, we prefer to single out some lemmas that we will employ in several situations.
For this reason, in this subsection we collect some results to be used later. 

Let $T>0$. 
We  introduce an auxiliary function defined by
\begin{equation}\label{eq:k}
k(t):=\left \{\begin{array}{l}
\displaystyle\sin \frac{\pi t}{T}\,\qquad\qquad \mbox{if}\,\, t\in\ [0,T]\,,\\
\\
0\,\qquad\qquad\quad\  \ \ \  \mbox{otherwise}\,.
\end{array}\right .
\end{equation}
In the following lemma we list some useful properties of $k$.

\begin{lemma} \label{th:k}
Set
\begin{equation}\label{eqn:K}
K(w):=\frac{T\pi}{\pi^2-T^2w^2}\,,\qquad w\in \C\,,
\end{equation}
the following properties hold.
\begin{itemize}
\item[(i)] 
 For any $w\in \C$ one has
 \begin{equation}\label{eqn:sinek2}
\overline{K(w)}=K(\overline{w})\,,
\qquad
\big|K(w)\big|=\big|K(\overline{w})\big|\,,
\end{equation}
\begin{equation}\label{eqn:sinek1}
\int_{0}^{\infty} k(t)e^{iw t}dt
= (1+e^{iw T})K(w)
\,.
\end{equation}
\item[(ii)] 
For any $z_i,w_i\in \C$, $i=1,2$, one has
\begin{multline}\label{eq:sinek2biss}
\int_{0}^{\infty} k(t)\Re(z_1e^{iw_1 t})\Re(z_2e^{iw_2t})dt
\\
=\frac12 \Re\Big(z_1z_2(1+e^{i(w_1+w_2) T})K(w_1+w_2)
+z_1\overline{z_2}(1+e^{i(w_1-\overline{w_2}) T})K(w_1-\overline{w_2})\Big)\,.
\end{multline}
\item[(iii)] 
Let $\overline{\gamma}>0$ and $j\in\N$. Then
for $T>2\pi/\overline{\gamma}$  and $w\in\C$, $|w|\ge\overline{\gamma} j$, one has
\begin{equation}\label{eq:sinek3}
\big|K(w)\big|\le 
\frac{4\pi}{T\overline{\gamma}^2(4j^2-1)}
\,.
\end{equation}
\end{itemize}

\end{lemma}
\begin{Proof}
(i) The proof is  straightforward. 

\noindent
(ii) We note that for any $z,w\in\C$ 
\begin{equation*}
\int_{0}^{\infty} k(t)\Re(ze^{iw t})dt
= \Re\big(z(1+e^{iw T})K(w)\big)\,.
\end{equation*}
Therefore, taking into account
\begin{equation*}
\Re(z_1e^{iw_1 t})\Re(z_2e^{iw_2t})=\frac12\Re\big(z_1z_2e^{i(w_1+w_2) t}+z_1\overline{z_2}e^{i(w_1-\overline{w_2}) t}\big)
\,,
\end{equation*}
it follows \eqref{eq:sinek2biss}.  

\noindent
(iii) We observe that
\begin{equation*}
\big|K(w)\big|=
\frac{\pi}{T\Big|w^2-\big(\frac{\pi}{T}\big)^2\Big|}
=\frac{4\pi}{T\overline{\gamma}^2\Big|4\big(\frac{w}{\overline{\gamma}}\big)^2-\big(\frac{2\pi}{T\overline{\gamma}}\big)^2\Big|}
\,.
\end{equation*}
Since $|w|\ge\overline{\gamma} j$ and $\frac{2\pi}{T\overline{\gamma}}<1$, 
we have
\begin{equation*}
\Big|4\Big(\frac{w}{\overline{\gamma}}\Big)^2-\Big(\frac{2\pi}{T\overline{\gamma}}\Big)^2\Big|
\ge
4\frac{|w|^2}{\overline{\gamma}^2}-\Big(\frac{2\pi}{T\overline{\gamma}}\Big)^2
\ge
4j^2-1\,,
\end{equation*}
and hence \eqref{eq:sinek3} holds true.
\end{Proof}
\begin{lemma}\label{le:gap}
If $\gamma>0$ is such that
\begin{equation*}
\liminf_{n\to\infty}\big(\Re\sigma_{n+1}-\Re\sigma_{n}\big)=\gamma\,,
\end{equation*}
then for any $\varepsilon\in (0,1)$ there exists $n_0\in\N$ such that
\begin{equation}\label{eq:Re_gap1}
|\Re\sigma_n-\Re\sigma_m|\ge \gamma\sqrt{1-\varepsilon} |n-m|\,,\qquad\forall n\,,m\ge n_0\,,
\end{equation}
\begin{equation}\label{eq:Re_gap2}
\Re\sigma_n\ge \gamma\sqrt{1-\varepsilon}\ n\,,\qquad\forall n\ge n_0\,.
\end{equation}
\end{lemma}

\begin{Proof}
For $\varepsilon\in (0,1)$ there exists $n_0\in\N$ such that
\begin{equation*}
\Re\sigma_{n+1}-\Re\sigma_n\ge\gamma\sqrt{1-\varepsilon}\qquad\quad\forall n\ge n_0\,,
\end{equation*}
whence \eqref{eq:Re_gap1} follows.
Moreover, in view of
\begin{equation}\label{eq:}
\liminf_{n\to\infty}\frac{\Re\sigma_{n+1}}{n+1}\ge
\liminf_{n\to\infty}\big(\Re\sigma_{n+1}-\Re\sigma_{n}\big)\,,
\end{equation}
see \cite[p. 54]{Ce}, \eqref{eq:Re_gap2} holds true.
\end{Proof}
\begin{lemma}\label{le:n0sum}
\begin{description}
\item[(i)]
For any $n_0\in\N$ and $n\ge n_0$ we have
 \begin{equation}\label{eq:telesc}
\sum_{\substack{m=n_0\\ m\not=n}}^\infty\
\frac{1}{4(m-n)^2-1}
\le
1\,.
\end{equation}
\item[(ii)]
Fixed $a,b\ge 0$ and $\varepsilon>0$, there exists $n_0\in\N$ large enough to satisfy
\begin{equation}\label{eq:n0sum}
\frac{a}{4n^2-1}+b\sum_{m=n_0}^\infty\frac{1}{4m^2-1}
\le\varepsilon
\qquad
\forall n\ge n_0
\,.
\end{equation}
\item[(iii)]
Fixed $a\ge 0$, $\nu>1/2$ and $\varepsilon>0$, there exists $n_0\in\N$ large enough to satisfy
\begin{equation}\label{eq:n0summ}
a\sum_{ n= n_0}^{\infty}\frac{1}{n^{2\nu}}
\le\varepsilon
\,.
\end{equation}
\end{description}
\end{lemma}
\begin{Proof}
(i) We have
 \begin{multline*}
\sum_{\substack{m=n_0\\ m\not=n}}^\infty\
\frac{1}{4(m-n)^2-1}
=
\sum_{m=n_0}^{n-1}\
\frac{1}{4(n-m)^2-1}
+\sum_{m=n+1}^\infty\
\frac{1}{4(m-n)^2-1}
\\
\le
2\sum_{j=1}^{\infty}\ \frac{1}{4j^2-1}
=\sum_{j=1}^{\infty}\
\Big(\frac{1}{2j-1}-\frac{1}{2j+1}\Big)=
1\,.
\end{multline*}

\noindent
(ii)
We observe that for $n\ge n_0$ we have
\begin{equation*}
4n^2-1
\ge 4n^{3/2}n_0^{1/2}-1
\ge n^{1/2}_0(4n^{3/2}-1)
\,,
\end{equation*}
and hence
\begin{equation*}
\frac{a}{4n^2-1}+b\sum_{m=n_0}^\infty\frac{1}{4m^2-1}
\le
\frac{1}{n^{1/2}_0}\bigg(a+b
\sum_{m=1}^\infty\frac{1}{4m^{3/2}-1}\bigg)
\,.
\end{equation*}
In conclusion, if one takes $n_0\in\N$ such that
\begin{equation*}
n_0
\ge
\frac{1}{\varepsilon^{2}}\bigg(a+b
\sum_{m=1}^\infty\frac{1}{4m^{3/2}-1}\bigg)^2
\,,
\end{equation*}
then \eqref{eq:n0sum} holds true.

\noindent
(iii) For $0<\delta<2\nu-1$ we have
\begin{equation*}
\sum_{ n= n_0}^{\infty}\frac{1}{n^{2\nu}}
\le
\frac{1}{n^{\delta}_0}\sum_{ n= 1}^{\infty}
\frac{1}{n^{2\nu-\delta}}\,,
\end{equation*}
whence, for $n_0\ge\bigg(\frac{a}\varepsilon\sum_{ n= 1}^{\infty}
\frac{1}{n^{2\nu-\delta}}\bigg)^{1/\delta}$ we have 
\eqref{eq:n0summ}.
\end{Proof}
\begin{lemma}\label{le:stimaK}
Suppose that
\begin{equation*}
\liminf_{n\to\infty}\big(\Re\sigma_{n+1}-\Re\sigma_{n}\big)=\gamma>0\,.
\end{equation*}
Then for any  $\varepsilon\in (0,1)$ and $T>\frac{2\pi}{\gamma\sqrt{1-\varepsilon}}$ there exists 
$n_0=n_0(\varepsilon)\in\N$ such that for any $n\ge n_0$ we have
\begin{equation}\label{eq:minus}
\sum_{\substack{m=n_0\\ m\not=n}}^\infty\ |K( \sigma_{n}-\overline{\sigma_m})|
+\sum_{m=n_0}^\infty|K( \sigma_n+\sigma_m)|
\le
\frac{4\pi}{T\gamma^2(1-\varepsilon)}\bigg(1+\sum_{m=n_0}^\infty\frac{1}{4m^{2}-1}\bigg)\,,
\end{equation}
\end{lemma}
\begin{Proof}
As regards the first inequality, we observe that, thanks to \eqref{eq:Re_gap1} and  \eqref{eq:sinek3}, for $\varepsilon\in (0,1)$ there exists $n_0\in\N$ such that
\begin{equation*}
\sum_{\substack{m=n_0\\ m\not=n}}^\infty\ |K( \sigma_{n}-\overline{\sigma_m})|
\le\frac{4\pi}{T\gamma^2(1-\varepsilon)}
\sum_{\substack{m=n_0\\ m\not=n}}^\infty\
\frac{1}{4(m-n)^2-1}
\,,
\end{equation*}
whence, in view of \eqref{eq:telesc} we get our statement.

Moreover, concerning the second estimate, thanks to \eqref{eq:Re_gap2}, we have
\begin{equation*}
| \sigma_n+\sigma_m|
\ge \Re\sigma_m\ge
\gamma\sqrt{1-\varepsilon}\ m\,,\qquad\forall m\ge n_0
\,.
\end{equation*}
Therefore, 
using again (\ref{eq:sinek3}) we obtain the required inequality.
\end{Proof}
The following result is an useful tool in the proof of the Ingham type inverse estimates. For the sake of completeness we prefer to give a detailed proof,  although it could be deduced from previous papers, see \cite{KL1}.

\begin{proposition}\label{pr:Fn}
Given any $\gamma>0$ suppose that
\begin{equation*}
\liminf_{n\to\infty}\big(\Re \sigma_{n+1}-\Re \sigma_{n}\big)=\gamma
\end{equation*}
and $\{F_n\}$ is a complex number sequence  such that 
$\sum_{n=1}^\infty\ |F_{n}|^2<+\infty$.

Then for any  $\varepsilon\in (0,1)$ and $T>\frac{2\pi}{\gamma\sqrt{1-\varepsilon}}$ 
 there exists 
$n_0=n_0(\varepsilon)\in\N$ independent of $T$ and $F_n$ such that we have
\begin{multline}\label{eq:Fn2'}
\int_{0}^{\infty} k(t)\Big| \sum_{n=n_0}^\infty  F_{n}e^{i \sigma_{n}t}+\overline{F_{n}}e^{-i \overline{\sigma_n}t}\Big|^2\ d t
\\
\ge
2\pi T\sum_{n=n_0}^\infty\bigg( \frac{1}{\pi^2+4T^2(\Im \sigma_{n})^2}
-
\frac{4}{T^2\gamma^2}(1+\varepsilon)\bigg) (1+e^{-2\Im \sigma_{n} T})|F_{n}|^2
\,,
\end{multline}
\begin{multline}\label{eq:Fn2''}
\int_{0}^{\infty} k(t)\Big| \sum_{n=n_0}^\infty  F_{n}e^{i \sigma_{n}t}+\overline{F_{n}}e^{-i \overline{\sigma_n}t}\Big|^2\ d t
\\
\le
2\pi T\sum_{n=n_0}^\infty\bigg( \frac{1}{\pi^2+4T^2(\Im \sigma_{n})^2}
+
\frac{4}{T^2\gamma^2}(1+\varepsilon)\bigg) (1+e^{-2\Im \sigma_{n} T})|F_{n}|^2
\,.
\end{multline}
\end{proposition}

\begin{Proof}
Let us first observe that 
\begin{equation*}
\Big| \sum_{n=n_0}^\infty F_{n}e^{i\sigma_{n}t}
+\overline{F_{n}}e^{-i \overline{\sigma_{n}}t}\Big|^2
=4\sum_{n, m=n_0}^\infty
\Re\big(F_{n}e^{i\sigma_{n}t}\big)\Re\big(F_{m}e^{i\sigma_{m}t}\big)
\,,
\end{equation*}
where $n_0\in\N$ will be chosen later.
From \eqref{eq:sinek2biss} we have
\begin{multline*}
\int_{0}^{\infty} k(t)\Big| 
\sum_{n=n_0}^\infty F_{n}e^{i\sigma_{n}t}
+\overline{F_{n}}e^{-i \overline{\sigma_{n}}t}\Big|^2\ d t
\\
=2\sum_{n, m=n_0}^\infty \Re \Big[F_{n}\overline{F_{m}} (1+e^{i( \sigma_{n}-\overline{\sigma_m}) T})
K( \sigma_{n}-\overline{\sigma_m})
+F_{n}F_{m} (1+e^{i( \sigma_{n}+\sigma_{m}) T})
K( \sigma_{n}+\sigma_{m})\Big]
\,.
\end{multline*}
Since \eqref{eqn:K} gives
$
\displaystyle
K( \sigma_{n}-\overline{\sigma_n})
=\frac{\pi T}{\pi^2+4T^2(\Im\sigma_{n})^2}\,,
$
it follows that
\begin{multline*}
\int_{0}^{\infty} k(t)\Big| \sum_{n=n_0}^\infty F_{n}e^{i 
\sigma_{n}t}+\overline{F_{n}}e^{-i \overline{\sigma_n}t}\Big|^2\ d t
-2\pi T\sum_{n=n_0}^\infty\ \frac{1+e^{-2\Im\sigma_{n} T}}{\pi^2+4T^2(\Im\sigma_{n})^2}|F_{n}|^2
\\
=2\sum_{\substack{n, m=n_0\\ n\not=m}}^\infty \Re \big[F_{n}\overline{F_{m}} (1+e^{i( \sigma_{n}-\overline{\sigma_m}) T})
K( \sigma_{n}-\overline{\sigma_m})\big]
+2\sum_{n, m=n_0}^\infty \Re \big[F_{n}F_{m} (1+e^{i( \sigma_{n}+\sigma_{m}) T})K( \sigma_{n}+\sigma_{m})\big]
\,.
\end{multline*}
Thus
\begin{multline}\label{eq:F2up}
\left|\int_{0}^{\infty} k(t)\Big| \sum_{n=n_0}^\infty F_{n}e^{i \sigma_{n}t}+\overline{F_{n}}e^{-i\overline{ \sigma_n}t}\Big|^2\ d t
-2\pi T\sum_{n=n_0}^\infty\ \frac{1+e^{-2\Im \sigma_{n} T}}{\pi^2+4T^2(\Im \sigma_{n})^2}|F_{n}|^2\right|
\\
\le
2\sum_{\substack{n, m=n_0\\ n\not=m}}^\infty |F_{n}| |F_{m}|
(1+e^{-\Im(\sigma_{n}+ \sigma_{m}) T})
|K( \sigma_{n}-\overline{\sigma_m})|
\\
+2\sum_{n, m=n_0}^\infty |F_{n}| |F_{m}|(1+e^{-\Im(\sigma_{n}+\sigma_{m}) T})
|K( \sigma_{n}+\sigma_{m})|
\,.
\end{multline}
By \eqref{eqn:sinek2} we have
$$
|K( \sigma_{n}-\overline{\sigma_m})|=|K( \sigma_{m}-\overline{\sigma_n})|\,,
$$
hence
\begin{multline*}
\sum_{\substack{n, m=n_0\\ n\not=m}}^\infty |F_{n}||F_{m}| |K( \sigma_{n}-\overline{\sigma_m})|
\le
\frac{1}{2}\sum_{\substack{n, m=n_0\\ n\not=m}}^\infty\ \big(|F_{n}|^2+ |F_{m}|^2\big) |K( \sigma_{n}-\overline{\sigma_m})|
\\
=
\frac{1}{2}\sum_{n=n_0}^\infty\ |F_{n}|^2\sum_{\substack{m=n_0\\ m\not=n}}^\infty\ |K( \sigma_{n}-\overline{\sigma_m})|
+\frac{1}{2}\sum_{m=n_0}^\infty\ |F_{m}|^2\sum_{\substack{n=n_0\\ n\not=m}}^\infty\ |K( \sigma_{m}-\overline{\sigma_n})|
\\
= \sum_{n=n_0}^\infty\ |F_{n}|^2\sum_{\substack{m=n_0\\ m\not=n}}^\infty\ |K( \sigma_{n}-\overline{\sigma_m})|
\,.
\end{multline*}
In the same manner we can see that 
\begin{equation*}
\sum_{\substack{n, m=n_0\\ n\not=m}}^\infty |F_{n}||F_{m}|e^{-\Im( \sigma_{n}+\sigma_{m}) T} |K( \sigma_{n}-\overline{\sigma_m})|
\le
 \sum_{n=n_0}^\infty\ e^{-2\Im\sigma_{n}T}|F_{n}|^2 \sum_{\substack{m=n_0\\ m\not=n}}^\infty\ |K( \sigma_{n}-\overline{\sigma_m})|
\,,
\end{equation*}
\begin{equation*}
\sum_{n, m=n_0}^\infty |F_{n}||F_{m}|(1+e^{-\Im( \sigma_{n}+\sigma_{m}) T}) |K( \sigma_{n}+\sigma_{m})|
\le
 \sum_{n=n_0}^\infty\ (1+e^{-2\Im \sigma_{n}T})|F_{n}|^2\sum_{m=n_0}^\infty\ |K( \sigma_{n}+\sigma_{m})|
\,.
\end{equation*}
Substituting these inequalities into  \eqref{eq:F2up} yields
\begin{multline*}
\left|\int_{0}^{\infty} k(t)\Big| \sum_{n=n_0}^\infty F_{n}e^{i \sigma_{n}t}+\overline{F_{n}}e^{-i \overline{\sigma_n}t}\Big|^2\ d t
-2\pi T\sum_{n=n_0}^\infty\ \frac{1+e^{-2\Im \sigma_{n} T}}{\pi^2+4T^2(\Im \sigma_{n})^2}|F_{n}|^2\right|
\\
\le
2\sum_{n=n_0}^\infty\ (1+e^{-2\Im \sigma_{n}T})|F_{n}|^2
\Bigg(\sum_{\substack{m=n_0\\ m\not=n}}^\infty\ |K( \sigma_{n}-\overline{\sigma_m})|
+\sum_{m=n_0}^\infty\ |K( \sigma_{n}+\sigma_{m})|\Bigg)
\,.
\end{multline*}
Fix now $\varepsilon\in (0,1)$ and $T>\frac{2\pi}{\gamma\sqrt{1-\varepsilon}}$. As for  $\varepsilon'\in (0,\varepsilon)$ one has $T>\frac{2\pi}{\gamma\sqrt{1-\varepsilon'}}$ too, we can employ Lemma \ref{le:stimaK} with $\varepsilon$ replaced by $\varepsilon'$. 
Thus taking $n_0$ as in Lemma \ref{le:stimaK} and applying \eqref{eq:minus}  we obtain
\begin{multline*}
\left|\int_{0}^{\infty} k(t)\Big| \sum_{n=n_0}^\infty F_{n}e^{i \sigma_{n}t}+\overline{F_{n}}e^{-i \overline{\sigma_n}t}\Big|^2\ d t
-2\pi T\sum_{n=n_0}^\infty\ \frac{1+e^{-2\Im\sigma_{n} T}}{\pi^2+4T^2(\Im\sigma_{n})^2}|F_{n}|^2\right|
\\
\le
\frac{8\pi}{T\gamma^2(1-\varepsilon')}\sum_{n=n_0}^\infty\ (1+e^{-2\Im\sigma_{n}T})|F_{n}|^2
\Bigg(1+\sum_{m=n_0}^\infty\frac{1}{4m^{2}-1}\Bigg)
\,.
\end{multline*}
By Lemma \ref{le:n0sum}-(ii) with $a=0$ and $b=1$ one can pick $n_0\in\N$ large enough to satisfy
\begin{equation*}
\sum_{m=n_0}^\infty\frac{1}{4m^{2}-1}
\le\varepsilon'
\,.
\end{equation*}
Therefore
\begin{multline*}
\left|\int_{0}^{\infty} k(t)\Big| \sum_{n=n_0}^\infty  F_{n}e^{i \sigma_{n}t}+\overline{F_{n}}e^{-i \overline{\sigma_n}t}\Big|^2\ d t
-2\pi T\sum_{n=n_0}^\infty\ \frac{1+e^{-2\Im\sigma_{n} T}}{\pi^2+4T^2(\Im\sigma_{n})^2}|F_{n}|^2\right|
\\
\le
\frac{8\pi}{T\gamma^2}\frac{1+\varepsilon'}{1-\varepsilon'}\sum_{n=n_0}^\infty\ (1+e^{-2\Im\sigma_{n} T})|F_{n}|^2
\,.
\end{multline*}
Taking $\varepsilon'\in (0,\varepsilon)$ such that $\frac{1+\varepsilon'}{1-\varepsilon'}<1+\varepsilon$, that is $\varepsilon'<\frac{\varepsilon}{2+\varepsilon}$, we obtain 
\begin{multline*}
\Bigg|\int_{0}^{\infty} k(t)\Big| \sum_{n=n_0}^\infty  F_{n}e^{i \sigma_{n}t}+\overline{F_{n}}e^{-i \overline{\sigma_n}t}\Big|^2\ d t
-2\pi T\sum_{n=n_0}^\infty\ \frac{1+e^{-2\Im\sigma_{n} T}}{\pi^2+4T^2(\Im\sigma_{n})^2}|F_{n}|^2\Bigg|
\\
\le
\frac{8\pi}{T\gamma^2}(1+\varepsilon)
\sum_{n=n_0}^\infty\ (1+e^{-2\Im\sigma_{n} T})|F_{n}|^2
\,,
\end{multline*}
which gives \eqref{eq:Fn2'} and \eqref{eq:Fn2''}.
\end{Proof}
\subsection{Inverse inequality}
Following the outline shown in Section \ref{se:outline} we have to estimate all three integrals on the right-hand side of \eqref{eq:f1+f2}. For this reason, for any term to bound we will establish
a corresponding lemma.
\begin{lemma}\label{le:}
For any  $\varepsilon\in (0,1)$ and $T>\frac{2\pi}{\gamma\sqrt{1-\varepsilon}}$ there exists 
$n_0=n_0(\varepsilon)\in\N$ independent of $T$ and $C_n$ such that 
we have
\begin{multline}\label{eq:u1-2u2}
\int_{0}^{\infty} k(t) \Big(\Big| \sum_{n=n_0}^{\infty}C_{n}e^{i\omega_{n} t}+\overline{C_{n}}e^{-i\overline{\omega_{n}}t}\Big|^2
-2\Big| \sum_{n=n_0}^\infty  c_nC_{n}e^{i \omega_{n}t}+\overline{c_nC_{n}}e^{-i \overline{\omega_n}t}\Big|^2\Big)\ d t
\\
\ge
2\pi T\sum_{n= n_0}^\infty
\Bigg(\frac{1-\varepsilon}{\pi^2+4T^2(\Im\omega_{n})^2}
-\frac{4}{T^2\gamma^2}(1+\varepsilon)
\Bigg)
(1+e^{-2\Im\omega_{n}T})|C_{n}|^2
\,.
\end{multline}
\end{lemma}

\begin{Proof}
Fix $\varepsilon\in (0,1)$ and $T>\frac{2\pi}{\gamma\sqrt{1-\varepsilon}}$.
Let us apply Proposition \ref{pr:Fn} with $\sigma_n=\omega_n$.
Indeed, for  $\varepsilon'\in (0,\varepsilon)$ to be chosen later  there exists 
$n_0$ independent of $T$ and $C_n$ such that 
from \eqref{eq:Fn2'} with $F_n=C_n$ and \eqref{eq:Fn2''} with $F_n=c_nC_n$ respectively we have 
\begin{multline}\label{eq:U1C}
\int_{0}^{\infty} k(t)\Big| \sum_{n=n_0}^\infty  C_{n}e^{i \omega_{n}t}+\overline{C_{n}}e^{-i \overline{\omega_n}t}\Big|^2\ d t
\\
\ge
2\pi T\sum_{n=n_0}^\infty\bigg( \frac{1}{\pi^2+4T^2(\Im \omega_{n})^2}
-
\frac{4}{T^2\gamma^2}(1+\varepsilon')\bigg) (1+e^{-2\Im \omega_{n} T})|C_{n}|^2
\,,
\end{multline}
\begin{multline}\label{eq:U2C}
\int_{0}^{\infty} k(t)\Big| \sum_{n=n_0}^\infty  c_nC_{n}e^{i \omega_{n}t}+\overline{c_nC_{n}}e^{-i \overline{\omega_n}t}\Big|^2\ d t
\\
\le
2\pi T\sum_{n=n_0}^\infty\bigg( \frac{1}{\pi^2+4T^2(\Im \omega_{n})^2}
+
\frac{4}{T^2\gamma^2}(1+\varepsilon')\bigg) (1+e^{-2\Im \omega_{n} T})|c_nC_{n}|^2
\,.
\end{multline}
Combining these inequalities gives
\begin{multline*}
\int_{0}^{\infty} k(t) \Big(\Big| \sum_{n=n_0}^{\infty}C_{n}e^{i\omega_{n} t}+\overline{C_{n}}e^{-i\overline{\omega_{n}}t}\Big|^2
-2\Big| \sum_{n=n_0}^\infty  c_nC_{n}e^{i \omega_{n}t}+\overline{c_nC_{n}}e^{-i \overline{\omega_n}t}\Big|^2\Big)\ d t
\\
\ge
2\pi T\sum_{n= n_0}^\infty
\Bigg(\frac{1-2|c_n|^2}{\pi^2+4T^2(\Im\omega_{n})^2}
-\frac{4}{T^2\gamma^2}(1+\varepsilon')
\big(1+2|c_n|^2\big)\Bigg)
(1+e^{-2\Im\omega_{n}T})|C_{n}|^2
\,.
\end{multline*}
We will choose $\varepsilon'$ in a suitable way to obtain our statement.
Thanks to \eqref{eq:cndn} for $n_0$ large enough we have $2|c_n|^2\le\varepsilon'$ for $n\ge n_0$. Hence
\begin{equation*}
(1+\varepsilon')\big(1+2|c_n|^2\big)\le(1+\varepsilon')^2\le1+3\varepsilon'
\hskip1cm
\forall n\ge n_0
\,.
\end{equation*}
Taking $\varepsilon'<\varepsilon/3$ yields
\begin{equation*}
(1+\varepsilon')\big(1+2|c_n|^2\big)\le1+\varepsilon
\hskip1cm
\forall n\ge n_0
\,.
\end{equation*}
Moreover, since
$2|c_n|^2\le\varepsilon$  we get
\eqref{eq:u1-2u2} and the proof is complete.
\end{Proof}
To estimate the second integral on the right-hand side of \eqref{eq:f1+f2} 
we state the following result, that  may be proved in much the same way 
as the previous lemma by means of Proposition \ref{pr:Fn} with 
$\sigma_n=\zeta_n$ and \eqref{eq:cndn}. For this reason we omit the proof.
\begin{lemma}\label{le:U2D}
For any  $\varepsilon\in (0,1)$ and $T>\frac{2\pi}{\gamma\sqrt{1-\varepsilon}}$ there exists 
$n_0=n_0(\varepsilon)\in\N$ independent of $T$ and $D_n$ such that 
we have
\begin{multline}\label{eq:U2D}
\int_{0}^{\infty} k(t) \Big(\Big| \sum_{n=n_0}^\infty d_nD_{n}e^{i \zeta_{n}t}+\overline{d_nD_{n}}e^{-i \overline{\zeta_n}t}\Big|^2
-4\Big| \sum_{n=n_0}^\infty D_{n}e^{i \zeta_{n}t}+\overline{D_{n}}e^{-i \overline{\zeta_n}t}\Big|^2\Big)\ d t
\\
\ge
2\pi T\sum_{n= n_0}^\infty
\Bigg(\frac{1-\varepsilon}{\pi^2+4T^2(\Im\zeta_{n})^2}
-\frac{4}{T^2\gamma^2}(1+\varepsilon)
\Bigg)
(1+e^{-2\Im\zeta_{n}T})|d_nD_{n}|^2
\,.
\end{multline}
\end{lemma}
Finally, we will give an estimate for the last integral on the right-hand side of \eqref{eq:f1+f2}.
\begin{lemma}
For any  $\varepsilon\in (0,1)$ and $T>0$ there exists 
$n_0=n_0(\varepsilon)\in\N$ independent of $T$ and $R_n$ such that 
we have
\begin{equation}\label{le:U1R0}
\int_{0}^{\infty} k(t) \Big|\sum_{n=n_0}^{\infty} R_{n} e^{r_n t}\Big|^2\ d t
\le
\varepsilon\ \pi T
\sum_{n= n_0}^\infty\frac{|C_{n}|^2+|d_nD_{n}|^2}{\pi^2+T^2r_{n}^2}
\,.
\end{equation}
\end{lemma}
\begin{Proof}
Our proof starts with the observation that
\eqref{eqn:sinek1} leads to
\begin{multline*}
\int_0^\infty k(t) \Big|\sum_{n=n_0}^{\infty} R_{n} e^{r_n t}\Big|^2\ dt=
\sum_{n,m=n_0}^{\infty} R_{n}R_{m} \int_0^\infty k(t) e^{(r_n+r_m) t}\ dt
\\
=
\sum_{n, m=n_0}^{\infty} R_{n} R_{m} (1+e^{(r_n+r_m )T})K(ir_n+ir_m)\,,
\end{multline*}
where $n_0\in\N$ has to be chosen later.
By the definition \eqref{eqn:K} of $K$  we have
\begin{equation*}
K(ir_n+ir_m)= \frac{T\pi}{\pi^2+T^2(r_n+r_m)^2}
\,.
\end{equation*}
Let us apply $r_n\le 0$ for $n\ge n'$ 
to  obtain
\begin{equation*}
1+e^{(r_n+r_m )T}\le2
\,.
\end{equation*}
Consequently, taking $n_0\ge n'$ we get 
\begin{equation*}
\int_0^\infty k(t) \Big|\sum_{n=n_0}^{\infty} R_{n} e^{r_n t}\Big|^2\ dt
\le
2\pi T\sum_{n, m= n_0}^{\infty} \frac{ |R_{n}| |R_{m}|}{\pi^2+T^2(r_n+r_m)^2}
\,.  
\end{equation*}
From
\eqref{eq:hom3}
we see that
\begin{multline*}
\int_0^\infty k(t) \Big|\sum_{n=n_0}^{\infty} R_{n} e^{r_n t}\Big|^2\ dt
\\
\le 
2\pi T\mu^2\sum_{n, m= n_0}^{\infty}\frac{\Big(|C_{n}|^2+ |d_nD_{n}|^2\Big)^{1/2}}{m^{\nu}}\ \frac{\Big(|C_{m}|^2+ |d_mD_{m}|^2\Big)^{1/2}}{n^{\nu}}
\frac{1}{\pi^2+T^2(r_n+r_m)^2}  
\,.  
\end{multline*}
Using again \eqref{eq:hom2} yields
\begin{multline*}
\sum_{n, m= n_0}^{\infty}\frac{\Big(|C_{n}|^2+ |d_nD_{n}|^2\Big)^{1/2}}{m^{\nu}}\ \frac{\Big(|C_{m}|^2+ |d_mD_{m}|^2\Big)^{1/2}}{n^{\nu}}
\frac{1}{\pi^2+T^2(r_n+r_m)^2}  
\\
\le
\frac12\sum_{ m= n_0}^{\infty}\frac{1}{m^{2\nu}}\sum_{ n= n_0}^{\infty}
\frac{|C_{n}|^2+ |d_nD_{n}|^2}{\pi^2+T^2r_n^2}
+\frac12\sum_{ n= n_0}^{\infty}\frac{1}{n^{2\nu}}\sum_{ m= n_0}^{\infty}
\frac{|C_{m}|^2+ |d_mD_{m}|^2}{\pi^2+T^2r_m^2}
\\
=
\sum_{ n= n_0}^{\infty}\frac{1}{n^{2\nu}}\sum_{ n= n_0}^{\infty}
\frac{|C_{n}|^2+ |d_nD_{n}|^2}{\pi^2+T^2r_n^2}
\,.  
\end{multline*}
Combining these inequalities we deduce that
\begin{equation*}
\int_0^\infty k(t) \Big|\sum_{n=n_0}^{\infty} R_{n} e^{r_n t}\Big|^2\ dt
\le
2\pi T\mu^2\sum_{ n= n_0}^{\infty}\frac{1}{n^{2\nu}}\sum_{ n= n_0}^{\infty}
\frac{|C_{n}|^2+ |d_nD_{n}|^2}{\pi^2+T^2r_n^2}\,.
\end{equation*}
Applying Lemma \ref{le:n0sum}-(iii) we conclude that
\eqref{le:U1R0} is proved.
\end{Proof}
We will establish the main result to obtain the inverse inequality. 
To simplify our notations, in the following we will use the  symbols
\begin{equation}\label{eq:notation}
\begin{split}
u_1^{n_0}(t) &:=\sum_{n=n_0}^{\infty}\Big(C_{n}e^{i\omega_{n} t}+\overline{C_{n}}e^{-i\overline{\omega_{n}}t}
+R_{n}e^{r_{n} t}+D_{n}e^{i\zeta_{n} t}+\overline{D_{n}}e^{-i\overline{\zeta_{n}}t}\Big)
\,,
\\
u_2^{n_0}(t) &:=\sum_{n=n_0}^{\infty}
\Big(d_nD_{n}e^{i\zeta_{n} t}+\overline{d_nD_{n}}e^{-i\overline{\zeta_{n}}t}
+c_nC_{n}e^{i\omega_{n} t}+\overline{c_nC_{n}}e^{-i\overline{\omega_{n}}t}\Big)
\,,
\end{split}
\end{equation}

\begin{theorem}\label{th:gamma>4alpha}
Assume $\gamma>4\alpha$.
Then for any $\varepsilon\in\big(0,\frac{\gamma^2-16\alpha^2}{\gamma^2+16\alpha^2}\big)$ and $T>\frac{2\pi}{\sqrt{\gamma^2(1-\varepsilon)-16\alpha^2(1+\varepsilon)}}$ there exist 
$n_0=n_0(\varepsilon)\in\N$, independent of $T$ and all coefficients of the series, and a constant $c(T,\varepsilon)>0$ such that 
\begin{multline}\label{eq:u1+u2}
\int_{0}^{\infty} k(t)\Big| \sum_{n=n_0}^{\infty}C_{n}e^{i\omega_{n} t}+\overline{C_{n}}e^{-i\overline{\omega_{n}}t}
+R_{n}e^{r_{n} t}+D_{n}e^{i\zeta_{n} t}+\overline{D_{n}}e^{-i\overline{\zeta_{n}}t}\Big|^2\ dt
\\
+\int_{0}^{\infty} k(t)\Big| \sum_{n=n_0}^{\infty}d_nD_{n}e^{i\zeta_{n} t}+\overline{d_nD_{n}}e^{-i\overline{\zeta_{n}}t}
+c_nC_{n}e^{i\omega_{n} t}+\overline{c_nC_{n}}e^{-i\overline{\omega_{n}}t}\Big|^2\ dt
\\
\ge
c(T,\varepsilon)\sum_{n= n_0}^\infty
(1+e^{-2\Im\omega_{n}T})\Big(|C_{n}|^2+ |d_nD_{n}|^2\Big)
\,.
\end{multline}
\end{theorem}
\begin{Proof}
Fix $\varepsilon\in (0,1)$, in view of \eqref{eq:notation}
our goal is to evaluate the following sum
\begin{equation}\label{eq:u1+u21}
\int_{0}^{\infty} k(t)
\big(|u_1^{n_0}(t)|^2+|u_2^{n_0}(t)|^2\big)\ dt
\,,
\end{equation}
where the index $n_0\in\N$ depending on  $\varepsilon$ will be chosen suitably. To this end, we bear in mind the comments given in Section \ref{se:outline}.
Indeed,
we observe that
\begin{multline*}
\int_{0}^{\infty} k(t)\Big| \sum_{n=n_0}^{\infty}C_{n}e^{i\omega_{n} t}+\overline{C_{n}}e^{-i\overline{\omega_{n}}t}
+R_{n}e^{r_{n} t}+D_{n}e^{i\zeta_{n} t}+\overline{D_{n}}e^{-i\overline{\zeta_{n}}t}\Big|^2\ dt
\\
\ge
\frac12\int_{0}^{\infty} k(t) \Big| \sum_{n=n_0}^{\infty}C_{n}e^{i\omega_{n} t}+\overline{C_{n}}e^{-i\overline{\omega_{n}}t}\Big|^2 \ d t
-2\int_{0}^{\infty} k(t)\Big| \sum_{n=n_0}^\infty D_{n}e^{i \zeta_{n}t}+\overline{D_{n}}e^{-i \overline{\zeta_n}t}\Big|^2\ d t
\\
-2\int_{0}^{\infty} k(t) \Big|\sum_{n=n_0}^{\infty} R_{n} e^{r_n t}\Big|^2\ d t
\end{multline*}
and
\begin{multline*}
\int_{0}^{\infty} k(t)\Big| \sum_{n=n_0}^{\infty}d_nD_{n}e^{i\zeta_{n} t}+\overline{d_nD_{n}}e^{-i\overline{\zeta_{n}}t}
+c_nC_{n}e^{i\omega_{n} t}+\overline{c_nC_{n}}e^{-i\overline{\omega_{n}}t}\Big|^2\ dt
\\
\ge
\frac12\int_{0}^{\infty} k(t) \Big| \sum_{n=n_0}^\infty d_nD_{n}e^{i \zeta_{n}t}+\overline{d_nD_{n}}e^{-i \overline{\zeta_n}t}\Big|^2\ dt
-\int_{0}^{\infty} k(t)\Big| \sum_{n=n_0}^\infty  c_nC_{n}e^{i \omega_{n}t}+\overline{c_nC_{n}}e^{-i \overline{\omega_n}t}\Big|^2\ d t
\,.
\end{multline*}
Combining these inequalities we obtain
\begin{multline*}
\int_{0}^{\infty} k(t)
\big(|u_1^{n_0}(t)|^2+|u_2^{n_0}(t)|^2\big)\ dt
\\
\ge
\frac12\int_{0}^{\infty} k(t) \Big(\Big| \sum_{n=n_0}^{\infty}C_{n}e^{i\omega_{n} t}+\overline{C_{n}}e^{-i\overline{\omega_{n}}t}\Big|^2
-2\Big| \sum_{n=n_0}^\infty  c_nC_{n}e^{i \omega_{n}t}+\overline{c_nC_{n}}e^{-i \overline{\omega_n}t}\Big|^2\Big)\ d t
\\
+\frac12\int_{0}^{\infty} k(t) \Big(\Big| \sum_{n=n_0}^\infty d_nD_{n}e^{i \zeta_{n}t}+\overline{d_nD_{n}}e^{-i \overline{\zeta_n}t}\Big|^2
-4\Big| \sum_{n=n_0}^\infty D_{n}e^{i \zeta_{n}t}+\overline{D_{n}}e^{-i \overline{\zeta_n}t}\Big|^2\Big)\ d t
\\
-2\int_{0}^{\infty} k(t) \Big|\sum_{n=n_0}^{\infty} R_{n} e^{r_n t}\Big|^2\ d t
\,.
\end{multline*}
We now take $T>\frac{2\pi}{\gamma\sqrt{1-\varepsilon}}$ to estimate the first two integrals on the right-hand side. We  introduce $\varepsilon'\in (0,\varepsilon)$ to choose suitably later. We also have $T>\frac{2\pi}{\gamma\sqrt{1-\varepsilon'}}$, so we can use \eqref{eq:u1-2u2} and \eqref{eq:U2D} respectively to obtain
\begin{multline*}
\int_{0}^{\infty} k(t) \Big(\Big| \sum_{n=n_0}^{\infty}C_{n}e^{i\omega_{n} t}+\overline{C_{n}}e^{-i\overline{\omega_{n}}t}\Big|^2
-2\Big| \sum_{n=n_0}^\infty  c_nC_{n}e^{i \omega_{n}t}+\overline{c_nC_{n}}e^{-i \overline{\omega_n}t}\Big|^2\Big)\ d t
\\
\ge
2\pi T\sum_{n= n_0}^\infty
\Bigg(\frac{1-\varepsilon'}{\pi^2+4T^2(\Im\omega_{n})^2}
-\frac{4}{T^2\gamma^2}(1+\varepsilon')
\Bigg)
(1+e^{-2\Im\omega_{n}T})|C_{n}|^2
\,,
\end{multline*}
\begin{multline*}
\int_{0}^{\infty} k(t) \Big(\Big| \sum_{n=n_0}^\infty d_nD_{n}e^{i \zeta_{n}t}+\overline{d_nD_{n}}e^{-i \overline{\zeta_n}t}\Big|^2
-4\Big| \sum_{n=n_0}^\infty D_{n}e^{i \zeta_{n}t}+\overline{D_{n}}e^{-i \overline{\zeta_n}t}\Big|^2\Big)\ d t
\\
\ge
2\pi T\sum_{n= n_0}^\infty
\Bigg(\frac{1-\varepsilon'}{\pi^2+4T^2(\Im\zeta_{n})^2}
-\frac{4}{T^2\gamma^2}(1+\varepsilon')
\Bigg)
(1+e^{-2\Im\zeta_{n}T})|d_nD_{n}|^2
\,.
\end{multline*}
By \eqref{eq:hom2} we get $|\Im\zeta_n|\le\Im\omega_n$ for $n\ge n_0$ with $n_0$ sufficiently large. Hence
\begin{equation*}
\frac{e^{-2\Im\zeta_{n}T}}{\pi^2+4T^2(\Im\zeta_{n})^2}
\ge
\frac{e^{-2\Im\omega_{n}T}}{\pi^2+4T^2(\Im\omega_{n})^2}
\qquad
\forall n\ge n_0
\,.
\end{equation*}
Therefore 
\begin{multline*}
\int_{0}^{\infty} k(t)
\big(|u_1^{n_0}(t)|^2+|u_2^{n_0}(t)|^2\big)\ dt
\\
\ge
\pi T\sum_{n= n_0}^\infty
\Bigg(\frac{1-\varepsilon'}{\pi^2+4T^2(\Im\omega_{n})^2}
-\frac{4}{T^2\gamma^2}(1+\varepsilon')
\Bigg)
(1+e^{-2\Im\omega_{n}T})\Big(|C_{n}|^2+ |d_nD_{n}|^2\Big)
\\
-2\int_{0}^{\infty} k(t) \Big|\sum_{n=n_0}^{\infty} R_{n} e^{r_n t}\Big|^2\ d t
\,.
\end{multline*}
Applying \eqref{le:U1R0} we obtain
\begin{multline}\label{eq:Irn}
\int_{0}^{\infty} k(t)
\big(|u_1^{n_0}(t)|^2+|u_2^{n_0}(t)|^2\big)\ dt
\\
\ge
\pi T\sum_{n= n_0}^\infty
\Bigg(\frac{1-\varepsilon'}{\pi^2+4T^2(\Im\omega_{n})^2}
-\frac{\varepsilon'}{\pi^2+T^2r_{n}^2}
-\frac{4}{T^2\gamma^2}(1+\varepsilon')
\Bigg)
(1+e^{-2\Im\omega_{n}T})\Big(|C_{n}|^2+ |d_nD_{n}|^2\Big)
\,.
\end{multline}
Now, we will choose $\varepsilon'\in (0,\varepsilon)$ such that for $n\ge n_0$ 
\begin{equation}\label{eq:TonT^2}
\frac{1-\varepsilon'}{\pi^2+4T^2(\Im\omega_{n})^2}
-\frac{\varepsilon'}{\pi^2+T^2r_{n}^2}
\ge
\frac{1-\varepsilon}{\pi^2+4T^2(\Im\omega_{n})^2}
\,,
\end{equation}
that is
\begin{equation*}
\frac{\varepsilon-\varepsilon'}{\pi^2+4T^2(\Im\omega_{n})^2}
-\frac{\varepsilon'}{\pi^2+T^2r_{n}^2}\ge0
\,,
\end{equation*}
\begin{equation*}
\pi^2(\varepsilon-2\varepsilon')
+T^2\big[(\varepsilon-\varepsilon')r_{n}^2
-4\varepsilon'(\Im\omega_{n})^2
\big]\ge0
\,.
\end{equation*}
To this end, we  need to have that
\begin{equation}\label{eq:epsi'}
\varepsilon-2\varepsilon'\ge0
\,,
\qquad
(\varepsilon-\varepsilon')r_{n}^2
-4\varepsilon'(\Im\omega_{n})^2
\ge0
\,.
\end{equation}
By \eqref{eq:hom2} for $n_0$ sufficiently large we have
\begin{equation*}
r_{n}^2\ge\frac{\chi^2}{2}\,,
\qquad
(\Im\omega_{n})^2\le\frac32\alpha^2\,.
\end{equation*}
Hence
\begin{equation*}
(\varepsilon-\varepsilon')r_{n}^2
-4\varepsilon'(\Im\omega_{n})^2
\ge
(\varepsilon-\varepsilon')\frac{\chi^2}2
-6\varepsilon'\alpha^2
\,.
\end{equation*}
Therefore taking
\begin{equation*}
\varepsilon'\le\min\Big\{\frac12,\frac{\chi^2}{\chi^2+12\alpha^2}\Big\}\varepsilon
\,,
\end{equation*}
we deduce \eqref{eq:epsi'}, and consequently \eqref{eq:TonT^2}.
So, from \eqref{eq:Irn} we have 
\begin{multline*}
\int_{0}^{\infty} k(t)
\big(|u_1^{n_0}(t)|^2+|u_2^{n_0}(t)|^2\big)\ dt
\\
\ge
\pi T\sum_{n= n_0}^\infty
\Bigg(\frac{1-\varepsilon}{\pi^2+4T^2(\Im\omega_{n})^2}
-\frac{4}{T^2\gamma^2}(1+\varepsilon)
\Bigg)
(1+e^{-2\Im\omega_{n}T})\Big(|C_{n}|^2+ |d_nD_{n}|^2\Big)
\,.
\end{multline*}
Since the previous inequality holds for any $\varepsilon\in(0,1)$, in particular it can be written for $\varepsilon'<\frac\varepsilon{2-\varepsilon}$, because 
this implies $\frac{1+\varepsilon'}{1-\varepsilon'}<\frac1{1-\varepsilon}$, and hence
\begin{equation*}
\frac{1-\varepsilon'}{\pi^2+4T^2(\Im\omega_{n})^2}
-\frac{4}{T^2\gamma^2}(1+\varepsilon')
\ge
(1-\varepsilon')\Bigg(\frac{1}{\pi^2+4T^2(\Im\omega_{n})^2}
-\frac{4}{T^2\gamma^2(1-\varepsilon)}
\Bigg)\,.
\end{equation*}
Therefore, taking also into account that $
(\Im\omega_{n})^2<\alpha^2(1+\varepsilon),
$
$n\ge n_0$, for $n_0$ large enough, we can write
\begin{multline}\label{eq:u1+u20}
\int_{0}^{\infty} k(t)
\big(|u_1^{n_0}(t)|^2+|u_2^{n_0}(t)|^2\big)\ dt
\\
\ge
\pi T(1-\varepsilon')\bigg(\frac{1}{\pi^2+4T^2\alpha^2(1+\varepsilon)}
-\frac{4}{T^2\gamma^2(1-\varepsilon)}\bigg)
\sum_{n= n_0}^\infty
(1+e^{-2\Im\omega_{n}T})\Big(|C_{n}|^2+ |d_nD_{n}|^2\Big)
\,.
\end{multline}
The constant 
\begin{equation*}
\frac{1}{\pi^2+4T^2\alpha^2(1+\varepsilon)}
-\frac{4}{T^2\gamma^2(1-\varepsilon)}
\end{equation*}
is positive if
\begin{equation}\label{eq:const}
T^2\big[\gamma^2(1-\varepsilon)-16\alpha^2(1+\varepsilon)\big]
>4\pi^2\,.
\end{equation}
Since $\gamma>4\alpha$ we have $\gamma^2(1-\varepsilon)-16\alpha^2(1+\varepsilon)>0$
if $\varepsilon<\frac{\gamma^2-16\alpha^2}{\gamma^2+16\alpha^2}$.
If we assume the more restrictive condition $T>\frac{2\pi}{\sqrt{\gamma^2(1-\varepsilon)-16\alpha^2(1+\varepsilon)}}$ with respect to that
$T>\frac{2\pi}{\gamma\sqrt{1-\varepsilon}}$, then \eqref{eq:const} holds true.
Finally, from \eqref{eq:u1+u20} and the definition \eqref{eq:u1+u21} of ${\cal I}_{n_0}$ we obtain \eqref{eq:u1+u2}.
%
%
\end{Proof}

We now observe that we can obtain a better estimate of the control time $T$
under an additional condition on the coefficients of the series.  Assuming  $|C_n|\le M |d_nD_{n}|$, we can follow the procedure sketched out at the end of Section \ref{se:outline} by using estimate \eqref{eq:f1+f2bis}. In particular,  to evaluate the term
$\int_{0}^{\infty}  k(t) | {\mathcal U}_{1}^C(t)|^2 d t$ we will employ the same trick used in \cite{LoretiSforza1}, giving first an estimate for $\int_{0}^{\infty} e^{2\alpha t} k(t) | {\mathcal U}_{1}^C(t)|^2 d t$ where $\displaystyle\alpha=\lim_{n\to\infty}{\Im}\om_n$ and then  multiplying by $e^{-2\alpha T}$ we will obtain the requested inequality.

\begin{theorem}\label{th:extracoe}
Assume 
\begin{equation}\label{eq:extracoe}
|C_n|\le M |d_nD_{n}|
\qquad\forall n\in\N
\,.
\end{equation}
Then, for any $\varepsilon\in (0,1)$ and $T>\frac{2\pi}{\gamma\sqrt{1-\varepsilon}}$ there exist 
$n_0=n_0(\varepsilon)\in\N$, independent of $T$ and all coefficients of the series, and a constant $c(T,\varepsilon)>0$ such that 
\begin{multline}\label{eq:u1+u2bis}
\int_{0}^{\infty} k(t)\Big| \sum_{n=n_0}^{\infty}C_{n}e^{i\omega_{n} t}+\overline{C_{n}}e^{-i\overline{\omega_{n}}t}
+R_{n}e^{r_{n} t}+D_{n}e^{i\zeta_{n} t}+\overline{D_{n}}e^{-i\overline{\zeta_{n}}t}\Big|^2\ dt
\\
+\int_{0}^{\infty} k(t)\Big| \sum_{n=n_0}^{\infty}d_nD_{n}e^{i\zeta_{n} t}+\overline{d_nD_{n}}e^{-i\overline{\zeta_{n}}t}
+c_nC_{n}e^{i\omega_{n} t}+\overline{c_nC_{n}}e^{-i\overline{\omega_{n}}t}\Big|^2\ dt
\\
\ge
c(T,\varepsilon)\sum_{n= n_0}^\infty
\Big(|C_{n}|^2+ |d_nD_{n}|^2\Big)
\,.
\end{multline}
\end{theorem}
\begin{Proof}
If $\displaystyle\alpha=\lim_{n\to\infty}{\Im}\om_n$, see \eqref{eq:hom2}, since
\begin{equation*}
\int_{0}^{\infty} e^{2\alpha t} k(t)\Big| \sum_{n=n_0}^\infty  C_{n}e^{i \omega_{n}t}+\overline{C_{n}}e^{-i \overline{\omega_n}t}\Big|^2\ d t
=\int_{0}^{\infty} k(t)\Big| \sum_{n=n_0}^\infty  C_{n}e^{i (\omega_{n}-i\alpha)t}+\overline{C_{n}}e^{-i \overline{(\omega_{n}-i\alpha)}t}\Big|^2\ d t
\,,
\end{equation*}
thanks to \eqref{eq:Fn2'}  we have
\begin{multline*}
\int_{0}^{\infty} e^{2\alpha t} k(t)\Big| \sum_{n=n_0}^\infty  C_{n}e^{i \omega_{n}t}+\overline{C_{n}}e^{-i \overline{\omega_n}t}\Big|^2\ d t
\\
\ge
2\pi T\sum_{n=n_0}^\infty\bigg( \frac{1}{\pi^2+4T^2(\Im \omega_{n}-\alpha)^2}
-\frac{4}{T^2\gamma^2}(1+\varepsilon')\bigg) (1+e^{-2(\Im \omega_{n}-\alpha) T})|C_{n}|^2
\,,
\end{multline*}
where $\varepsilon'\in(0,\varepsilon)$ will be chosen later.
Therefore, multiplying by $e^{-2\alpha T}$ and taking into account the definition \eqref{eq:k} of the function $k$, we get
\begin{multline*}
\int_{0}^{\infty} k(t)\Big| \sum_{n=n_0}^\infty  C_{n}e^{i \omega_{n}t}+\overline{C_{n}}e^{-i \overline{\omega_n}t}\Big|^2\ d t
\\
\ge
2\pi Te^{-2\alpha T}\sum_{n=n_0}^\infty\bigg( \frac{1}{\pi^2+4T^2(\Im \omega_{n}-\alpha)^2}
-\frac{4}{T^2\gamma^2}(1+\varepsilon')\bigg) (1+e^{-2(\Im \omega_{n}-\alpha) T})|C_{n}|^2
\,.
\end{multline*}
Now, we can take $4(\Im \omega_{n}-\alpha)^2<\gamma^2\varepsilon/8$ for $n\ge n_0$ and $1+\varepsilon'<\frac1{1-\varepsilon/2}$ for $\varepsilon'<\frac\varepsilon{2-\varepsilon}$, to have
\begin{multline*}
\frac{1}{\pi^2+4T^2(\Im \omega_{n}-\alpha)^2}-\frac{4}{T^2\gamma^2}(1+\varepsilon')
\\
>
\frac{1}{\pi^2+T^2\gamma^2\varepsilon/8}-\frac{4}{T^2\gamma^2(1-\varepsilon/2)}
=\frac{T^2\gamma^2(1-\varepsilon)-4\pi^2}{(\pi^2+T^2\gamma^2\varepsilon/8)T^2\gamma^2(1-\varepsilon/2)}
\end{multline*}
and $T^2\gamma^2(1-\varepsilon)-4\pi^2>0$ for $T>\frac{2\pi}{\gamma\sqrt{1-\varepsilon}}$.
So, we get
\begin{multline}\label{eq:onlyC1}
\int_{0}^{\infty} k(t)\Big| \sum_{n=n_0}^\infty  C_{n}e^{i \omega_{n}t}+\overline{C_{n}}e^{-i \overline{\omega_n}t}\Big|^2\ d t
\\
\ge
2\pi Te^{-2\alpha T}\frac{T^2\gamma^2(1-\varepsilon)-4\pi^2}{(\pi^2+T^2\gamma^2\varepsilon/8)T^2\gamma^2(1-\varepsilon/2)}\sum_{n=n_0}^\infty (1+e^{-2(\Im \omega_{n}-\alpha) T})|C_{n}|^2
\,.
\end{multline}
On the other hand,  from \eqref{eq:Fn2''} it follows
\begin{multline*}
\int_{0}^{\infty} k(t)\Big| \sum_{n=n_0}^\infty  c_nC_{n}e^{i \omega_{n}t}+\overline{c_nC_{n}}e^{-i \overline{\omega_n}t}\Big|^2\ d t
\\
\le
2\pi T\sum_{n=n_0}^\infty\bigg( \frac{1}{\pi^2+4T^2(\Im \omega_{n})^2}
+
\frac{4}{T^2\gamma^2}(1+\varepsilon')\bigg) (1+e^{-2\Im \omega_{n} T})|c_nC_{n}|^2
\\
\le
2\pi T\sum_{n=n_0}^\infty  M|c_n|^2\bigg( \frac{1}{\pi^2+4T^2(\Im \zeta_{n})^2}
+
\frac{4}{T^2\gamma^2}(1+\varepsilon')\bigg) (1+e^{-2\Im \zeta_{n} T})|d_nD_{n}|^2
\,,
\end{multline*}
thanks also to $\Im\omega_n\ge|\Im\zeta_n|$ and $|C_n|\le M |d_nD_{n}|$.
Moreover, again by \eqref{eq:Fn2''} and  the previous inequality we have
\begin{multline*}
\int_{0}^{\infty} k(t) \Big(
2\Big| \sum_{n=n_0}^\infty D_{n}e^{i \zeta_{n}t}+\overline{D_{n}}e^{-i \overline{\zeta_n}t}\Big|^2
+\Big| \sum_{n=n_0}^\infty  c_nC_{n}e^{i \omega_{n}t}+\overline{c_nC_{n}}e^{-i \overline{\omega_n}t}\Big|^2\Big)
\ d t
\\
\le
2\pi T\sum_{n=n_0}^\infty  \Big(\frac2{|d_n|^2}+M|c_n|^2\Big)\bigg( \frac{1}{\pi^2+4T^2(\Im \zeta_{n})^2}
+
\frac{4}{T^2\gamma^2}(1+\varepsilon')\bigg) (1+e^{-2\Im \zeta_{n} T})|d_nD_{n}|^2
\,.
\end{multline*}
Choosing $n_0$ sufficiently large such that $\frac2{|d_n|^2}+M|c_n|^2\le\varepsilon'$ for any $n\ge n_0$, from the above estimate we deduce
\begin{multline}\label{eq:U1D2C}
\int_{0}^{\infty} k(t) \Big(
2\Big| \sum_{n=n_0}^\infty D_{n}e^{i \zeta_{n}t}+\overline{D_{n}}e^{-i \overline{\zeta_n}t}\Big|^2
+\Big| \sum_{n=n_0}^\infty  c_nC_{n}e^{i \omega_{n}t}+\overline{c_nC_{n}}e^{-i \overline{\omega_n}t}\Big|^2\Big)
\ d t
\\
\le
2\pi T\varepsilon'\sum_{n=n_0}^\infty  \bigg( \frac{1}{\pi^2+4T^2(\Im \zeta_{n})^2}
+
\frac{4}{T^2\gamma^2}(1+\varepsilon')\bigg) (1+e^{-2\Im \zeta_{n} T})|d_nD_{n}|^2
\,.
\end{multline}
In addition, from \eqref{le:U1R0}, using again  $|C_n|\le M |d_nD_{n}|$  and \eqref{eq:hom2} we get
\begin{equation}\label{eq:U1R1}
\int_{0}^{\infty} k(t) \Big|\sum_{n=n_0}^{\infty} R_{n} e^{r_n t}\Big|^2\ d t
\le
 \pi T \varepsilon'
\sum_{n= n_0}^\infty\frac{|d_nD_{n}|^2}{\pi^2+T^2r_{n}^2}
\le
 \pi T \varepsilon'
\sum_{n= n_0}^\infty\frac{|d_nD_{n}|^2}{\pi^2+4T^2(\Im\zeta_{n})^2}
\,.
\end{equation}
Combining \eqref{eq:U1D2C} and \eqref{eq:U1R1} (with $\varepsilon'$ replaced by $\varepsilon'/2$) we obtain
\begin{multline}\label{eq:D1C2R}
\int_{0}^{\infty} k(t) \Big(
2\Big| \sum_{n=n_0}^\infty D_{n}e^{i \zeta_{n}t}+\overline{D_{n}}e^{-i \overline{\zeta_n}t}\Big|^2
+\Big| \sum_{n=n_0}^\infty  c_nC_{n}e^{i \omega_{n}t}+\overline{c_nC_{n}}e^{-i \overline{\omega_n}t}\Big|^2
+2 \Big|\sum_{n=n_0}^{\infty} R_{n} e^{r_n t}\Big|^2\Big)
\ d t
\\
\le
2\pi T\varepsilon'\sum_{n=n_0}^\infty  \bigg( \frac{1}{\pi^2+4T^2(\Im \zeta_{n})^2}
+
\frac{4}{T^2\gamma^2}(1+\varepsilon')\bigg) (1+e^{-2\Im \zeta_{n} T})|d_nD_{n}|^2
\,.
\end{multline}
In virtue of \eqref{eq:Fn2'}  we get
\begin{multline*}
\int_{0}^{\infty} k(t) \Big| \sum_{n=n_0}^\infty d_nD_{n}e^{i \zeta_{n}t}+\overline{d_nD_{n}}e^{-i \overline{\zeta_n}t}\Big|^2
\ d t
\\
\ge
2\pi T\sum_{n= n_0}^\infty
\Bigg(\frac{1}{\pi^2+4T^2(\Im\zeta_{n})^2}
-\frac{4}{T^2\gamma^2}(1+\varepsilon')
\Bigg)
(1+e^{-2\Im\zeta_{n}T})|d_nD_{n}|^2
\,.
\end{multline*}
From the above formula and \eqref{eq:D1C2R}, taking  $\varepsilon'\le\varepsilon/3$ but writing again $\varepsilon'$ instead of $\varepsilon$, we have
\begin{multline*}
\int_{0}^{\infty} k(t) \Big| \sum_{n=n_0}^\infty d_nD_{n}e^{i \zeta_{n}t}+\overline{d_nD_{n}}e^{-i \overline{\zeta_n}t}\Big|^2
\ d t
\\
-2
\int_{0}^{\infty} k(t) \Big(
2\Big| \sum_{n=n_0}^\infty D_{n}e^{i \zeta_{n}t}+\overline{D_{n}}e^{-i \overline{\zeta_n}t}\Big|^2
+\Big| \sum_{n=n_0}^\infty  c_nC_{n}e^{i \omega_{n}t}+\overline{c_nC_{n}}e^{-i \overline{\omega_n}t}\Big|^2
+2 \Big|\sum_{n=n_0}^{\infty} R_{n} e^{r_n t}\Big|^2\Big)
\ d t
\\
\ge
2\pi T\sum_{n= n_0}^\infty
\Bigg(\frac{1-\varepsilon'}{\pi^2+4T^2(\Im\zeta_{n})^2}
-\frac{4}{T^2\gamma^2}(1+\varepsilon')
\Bigg)
(1+e^{-2\Im\zeta_{n}T})|d_nD_{n}|^2
\,.
\end{multline*}
Taking $4(\Im\zeta_{n})^2<\gamma^2\varepsilon/8$ for $n\ge n_0$ and $\frac{1+\varepsilon'}{1-\varepsilon'}<\frac1{1-\varepsilon/2}$ for $\varepsilon'<\frac\varepsilon{4-\varepsilon}$ yields
\begin{multline*}
\frac{1-\varepsilon'}{\pi^2+4T^2(\Im\zeta_{n})^2}
-\frac{4}{T^2\gamma^2}(1+\varepsilon')
=(1-\varepsilon')
\bigg(\frac{1}{\pi^2+4T^2(\Im\zeta_{n})^2}
-\frac{4(1+\varepsilon')}{T^2\gamma^2(1-\varepsilon')}
\bigg)
\\
\ge
(1-\varepsilon')
\bigg(\frac{1}{\pi^2+T^2\gamma^2\varepsilon/8}
-\frac{4}{T^2\gamma^2(1-\varepsilon/2)}
\bigg)
=
(1-\varepsilon')
\bigg(\frac{T^2\gamma^2(1-\varepsilon)-4\pi^2}{(\pi^2+T^2\gamma^2\varepsilon/8)T^2\gamma^2(1-\varepsilon/2)}
\bigg)
\,.
\end{multline*}
Therefore,  for $T>\frac{2\pi}{\gamma\sqrt{1-\varepsilon}}$
we obtain
\begin{multline*}
\int_{0}^{\infty} k(t) \Big| \sum_{n=n_0}^\infty d_nD_{n}e^{i \zeta_{n}t}+\overline{d_nD_{n}}e^{-i \overline{\zeta_n}t}\Big|^2
\ d t
\\
-2
\int_{0}^{\infty} k(t) \Big(
2\Big| \sum_{n=n_0}^\infty D_{n}e^{i \zeta_{n}t}+\overline{D_{n}}e^{-i \overline{\zeta_n}t}\Big|^2
+\Big| \sum_{n=n_0}^\infty  c_nC_{n}e^{i \omega_{n}t}+\overline{c_nC_{n}}e^{-i \overline{\omega_n}t}\Big|^2
+2 \Big|\sum_{n=n_0}^{\infty} R_{n} e^{r_n t}\Big|^2\Big)
\ d t
\\
\ge
2\pi T(1-\varepsilon)
\bigg(\frac{T^2\gamma^2(1-\varepsilon)-4\pi^2}{(\pi^2+T^2\gamma^2\varepsilon/8)T^2\gamma^2(1-\varepsilon/2)}
\bigg)
\sum_{n= n_0}^\infty
(1+e^{-2\Im\zeta_{n}T})|d_nD_{n}|^2
\,.
\end{multline*}
In conclusion, for any $T>\frac{2\pi}{\gamma\sqrt{1-\varepsilon}}$,
combining the previous estimate with \eqref{eq:onlyC1} gives 
\begin{multline*}
\int_{0}^{\infty} k(t)
\big(|u_1^{n_0}(t)|^2+|u_2^{n_0}(t)|^2\big)\ dt
\\
\ge
\pi T\min\{e^{-2\alpha T},(1-\varepsilon)\}
\bigg(\frac{T^2\gamma^2(1-\varepsilon)-4\pi^2}{(\pi^2+T^2\gamma^2\varepsilon/8)T^2\gamma^2(1-\varepsilon/2)}
\bigg)
\sum_{n=n_0}^\infty\Big(|C_{n}|^2 +|d_nD_{n}|^2\Big)
\,,
\end{multline*}
that is \eqref{eq:u1+u2bis}.
\end{Proof}

%
%

\subsection{Direct inequality}
As for the inverse inequality, to prove  direct estimates we need to introduce an auxiliary function. Let $T>0$ and define  
\begin{equation}\label{eq:kcos}
k^*(t):=\left \{\begin{array}{l}
\cos \frac{\pi t}{2T}\,\qquad\qquad \mbox{if}\ |t|\le T\,,\\
\\
0\,\qquad\qquad\quad\  \ \ \  \mbox{if}\  |t|>T\,.
\end{array}\right .
\end{equation}
For the sake of completeness, we list some standard properties of $k^*$ in the following lemma.
\begin{lemma} \label{th:k}
Set
\begin{equation}\label{eqn:K*}
K^*(u):=\frac{4T\pi}{\pi^2-4T^2u^2}\,,\qquad u\in \C\,,
\end{equation}
the following properties hold  for any $u\in \C$
\begin{equation}\label{eqn:k1}
\int_{-\infty}^{\infty} k^*(t)e^{iu t}dt=\cos(uT)K^*(u)\,,
\end{equation}
\begin{equation}\label{eqn:k2bis}
\overline{K^*(u)}=K^*(\overline{u})\,,
\quad
\big|K^*(u)\big|=\big|K^*(\overline{u})\big|.
\end{equation}
Set 
$K_{T}(u)=\frac{T\pi}{\pi^2-T^2u^2}$ we have
\begin{equation}\label{eqn:k3}
K^*(u)=2K_{2T}(u)\,.
\end{equation}
Moreover for any $z_i,w_i\in \C$, $i=1,2$, one has
\begin{multline}\label{eq:cosbiss}
\int_{-\infty}^{\infty} k^*(t)\Re(z_1e^{iw_1 t})\Re(z_2e^{iw_2t})dt
\\
=\frac12 \Re\Big(z_1z_2\cos((w_1+w_2) T)K(w_1+w_2)
+z_1\overline{z_2}\cos((w_1-\overline{w_2}) T)K(w_1-\overline{w_2})\Big)\,.
\end{multline}\end{lemma}
From now on we will denote with $c(T)$ a positive constant depending on $T$. 
\begin{proposition}\label{pr:Fnd}
Let $\gamma>0$. Suppose that $\{\sigma_{n}\}$ is a complex number sequence satisfying
\begin{equation*}
\liminf_{n\to\infty}\big(\Re \sigma_{n+1}-\Re \sigma_{n}\big)=\gamma\,,
\qquad
\{{\Im}\sigma_n\}
\quad
\mbox{bounded}.
\end{equation*}
Then for any complex number sequence $\{F_n\}$ with 
$\sum_{n=1}^\infty\ |F_{n}|^2<+\infty$,
$\varepsilon\in (0,1)$ and $T>\frac{\pi}{\gamma\sqrt{1-\varepsilon}}$ 
 there exist $c(T)>0$ and
$n_0=n_0(\varepsilon)\in\N$ independent of $T$ and $F_n$ such that 
\begin{equation}\label{eq:Fndir}
\int_{-\infty}^{\infty} k^*(t)\Big| \sum_{n=n_0}^\infty  F_{n}e^{i \sigma_{n}t}+\overline{F_{n}}e^{-i \overline{\sigma_n}t}\Big|^2\ d t
\le
c(T)\sum_{n=n_0}^\infty|F_{n}|^2
\,.
\end{equation}
\end{proposition}

\begin{Proof}
Let us first observe that 
\begin{equation*}
\Big| \sum_{n=n_0}^\infty F_{n}e^{i\sigma_{n}t}
+\overline{F_{n}}e^{-i \overline{\sigma_{n}}t}\Big|^2
=4\sum_{n, m=n_0}^\infty
\Re\big(F_{n}e^{i\sigma_{n}t}\big)\Re\big(F_{m}e^{i\sigma_{m}t}\big)
\,,
\end{equation*}
where the index $n_0\in\N$ depending on  $\varepsilon$ will be chosen later.
From \eqref{eq:cosbiss} we have
\begin{multline*}
\int_{-\infty}^{\infty} k^*(t)\Big| 
\sum_{n=n_0}^\infty F_{n}e^{i\sigma_{n}t}
+\overline{F_{n}}e^{-i \overline{\sigma_{n}}t}\Big|^2\ d t
\\
=2\sum_{n, m=n_0}^\infty \Re \Big[F_{n}\overline{F_{m}} \cos(( \sigma_{n}-\overline{\sigma_m}) T)
K^*( \sigma_{n}-\overline{\sigma_m})
+F_{n}F_{m} \cos(( \sigma_{n}+\sigma_{m}) T)
K^*( \sigma_{n}+\sigma_{m})\Big]
\,.
\end{multline*}
Applying  the elementary estimates $\Re z\le|z|$ and
$|\cos z|\le \cosh (\Im z)$, $z\in\C$,
we obtain
\begin{multline*}
\int_{-\infty}^{\infty} k^*(t)\Big| \sum_{n=n_0}^\infty F_{n}e^{i \sigma_{n}t}+\overline{F_{n}}e^{-i\overline{ \sigma_n}t}\Big|^2\ d t
\\
\le
2\sum_{n, m=n_0}^\infty |F_{n}| |F_{m}|
\cosh(\Im(\sigma_{n}+ \sigma_{m}) T)
\big[|K^*( \sigma_{n}-\overline{\sigma_m})|
+|K^*( \sigma_{n}+\sigma_{m})|\big]
\,.
\end{multline*}
Since the sequence $\{{\Im}\sigma_n\}$ is bounded
we have
\begin{equation*}
\cosh(\Im(\sigma_{n}+ \sigma_{m}) T)\le e^{2T\sup|\Im\sigma_n|}
\qquad
\forall  n,m\in\N\,.
\end{equation*}
Hence
\begin{multline*}
\int_{-\infty}^{\infty} k^*(t)\Big| \sum_{n=n_0}^\infty F_{n}e^{i \sigma_{n}t}+\overline{F_{n}}e^{-i\overline{ \sigma_n}t}\Big|^2\ d t
\\
\le
2e^{2T\sup|\Im\sigma_n|}\sum_{n, m=n_0}^\infty |F_{n}| |F_{m}|
\big[|K^*( \sigma_{n}-\overline{\sigma_m})|
+|K^*( \sigma_{n}+\sigma_{m})|\big]
\,.
\end{multline*}
Thanks to \eqref{eqn:k2bis} we get
$
|K^*( \sigma_n-\overline{\sigma_m})|=|K^*( \sigma_m-\overline{\sigma_n})|\,.
$
Therefore
\begin{multline*}
\int_{-\infty}^{\infty} k^*(t)\Big| \sum_{n=n_0}^\infty F_{n}e^{i \sigma_{n}t}+\overline{F_{n}}e^{-i\overline{ \sigma_n}t}\Big|^2\ d t
\\
\le
2e^{2T\sup|\Im\sigma_n|}\sum_{n=n_0}^\infty |F_{n}|^2
\sum_{ m=n_0}^\infty \big[|K^*( \sigma_{n}-\overline{\sigma_m})|
+|K^*( \sigma_{n}+\sigma_{m})|\big]
\,.
\end{multline*}
Since \eqref{eqn:K*} gives
\begin{equation*}
K^*( \sigma_{n}-\overline{\sigma_n})
=\frac{4\pi T}{\pi^2+16T^2(\Im\sigma_{n})^2}
\le\frac{4 T}{\pi}
\,,
\end{equation*}
it follows that
\begin{multline}\label{eq:Fndir0}
\int_{-\infty}^{\infty} k(t)\Big| \sum_{n=n_0}^\infty F_{n}e^{i 
\sigma_{n}t}+\overline{F_{n}}e^{-i \overline{\sigma_n}t}\Big|^2\ d t
\le
\frac8{\pi} e^{2T\sup|\Im\sigma_n|}T\sum_{n=n_0}^\infty\ |F_{n}|^2
\\
+2e^{2T\sup|\Im\sigma_n|}\sum_{n=n_0}^\infty |F_{n}|^2
\Big[\sum_{\substack{m=n_0\\ m\not=n}}^\infty |K^*( \sigma_{n}-\overline{\sigma_m})|
+\sum_{m=n_0}^\infty K^*( \sigma_{n}+\sigma_{m})\Big]
\,.
\end{multline}
Note that by \eqref{eqn:k3} we can apply Lemma \ref{le:stimaK}: for any
$\varepsilon\in (0,1)$ and $2T>\frac{2\pi}{\gamma\sqrt{1-\varepsilon}}$ there exists $n_0\in\N$ such that 
\begin{equation*}
\sum_{\substack{m=n_0\\ m\not=n}}^\infty |K^*( \sigma_{n}-\overline{\sigma_m})|
+\sum_{m=n_0}^\infty K^*( \sigma_{n}+\sigma_{m})
\le
\frac{{2\pi}}{T\gamma^2(1-\varepsilon)}
\Big(1+\sum_{n=1}^\infty\frac{1}{4n^{2}-1}\Big)
\,.
\end{equation*}
Substituting the previous estimate into \eqref{eq:Fndir0} gives \eqref{eq:Fndir}.
\end{Proof}
\begin{proposition}\label{pr:Rnd}
For any  $n_0\in\N$,  $n_0\ge n'$, and $T>0$ there exists $c(T)>0$ such that 
\begin{equation}\label{le:U1R}
\int_{-\infty}^\infty k^*(t) \bigg|\sum_{n=n_0}^{\infty} R_{n} e^{r_n t}\bigg|^2\ dt
\le
c(T)
\sum_{n= n_0}^\infty
\Big(|C_{n}|^2+ |d_nD_{n}|^2\Big)
\,.
\end{equation}
\end{proposition}
\begin{Proof}
Fixed  $n_0\in\N$,  $n_0\ge n'$, we observe that
\eqref{eqn:k1} leads to
\begin{multline*}
\int_{-\infty}^\infty k^*(t) \bigg|\sum_{n=n_0}^{\infty} R_{n} e^{r_n t}\bigg|^2\ dt=
\sum_{n,m=n_0}^{\infty} R_{n}R_{m} \int_{\infty}^\infty k^*(t) e^{(r_n+r_m) t}\ dt
\\
=
\sum_{n, m=n_0}^{\infty} R_{n} R_{m} \cosh((r_n+r_m )T)K^*(ir_n+ir_m)\,.
\end{multline*}
By the definition \eqref{eqn:K*} of $K^*$  we have
\begin{equation*}
K^*(ir_n+ir_m)= \frac{4\pi T}{\pi^2+4T^2(r_n+r_m)^2}
\le
\frac{4 T}{\pi}
\,.
\end{equation*}
In addition, since the sequence $\{r_n\}$ is bounded
we have
\begin{equation*}
\cosh((r_{n}+ r_{m}) T)\le e^{2T\sup|r_n|}
\qquad
\forall  n,m\in\N\,.
\end{equation*}
Consequently, 
\begin{equation*}
\int_{-\infty}^\infty k^*(t) \bigg|\sum_{n=n_0}^{\infty} R_{n} e^{r_n t}\bigg|^2\ dt
\le
\frac{4 T}{\pi} e^{2T\sup|r_n|}\sum_{n, m=n_0}^{\infty} |R_{n}| |R_{m}|
\,.  
\end{equation*}
Since $n_0\ge n'$, by
\eqref{eq:hom3} 
we have that
\begin{equation*}
\int_{-\infty}^\infty k^*(t) \bigg|\sum_{n=n_0}^{\infty} R_{n} e^{r_n t}\bigg|^2\ dt
\le 
\frac{4 T}{\pi} e^{2T\sup|r_n|}\sum_{n, m= n_0}^{\infty}\frac{\Big(|C_{n}|^2+ |d_nD_{n}|^2\Big)^{1/2}}{m^{\nu}}\ \frac{\Big(|C_{m}|^2+ |d_mD_{m}|^2\Big)^{1/2}}{n^{\nu}}
\,.  
\end{equation*}
Moreover
\begin{multline*}
\sum_{n, m= n_0}^{\infty}
\frac{\Big(|C_{n}|^2+ |d_nD_{n}|^2\Big)^{1/2}}{m^{\nu}}\ \frac{\Big(|C_{m}|^2+ |d_mD_{m}|^2\Big)^{1/2}}{n^{\nu}}
\\
\le
\frac12\sum_{ m= n_0}^{\infty}\frac{1}{m^{2\nu}}
\sum_{ n= n_0}^{\infty}
\Big(|C_{n}|^2+ |d_nD_{n}|^2\Big)
+\frac12\sum_{ n= n_0}^{\infty}\frac{1}{n^{2\nu}}
\sum_{ m= n_0}^{\infty}
\Big(|C_{n}|^2+ |d_nD_{n}|^2\Big)
\\
=
\sum_{ n= 1}^{\infty}\frac{1}{n^{2\nu}}
\sum_{ n= n_0}^{\infty}
\Big(|C_{n}|^2+ |d_nD_{n}|^2\Big)
\,.  
\end{multline*}
Combining these inequalities  we conclude that
\eqref{le:U1R} is proved.
\end{Proof}
\begin{theorem}\label{th:Diretta}
For any $\varepsilon\in (0,1)$ and 
$T>\frac{\pi}{\ga\sqrt{1-\varepsilon}}$ there exist $n_0=n_0(\varepsilon)\in\N$ and $c(T)>0$ such that 
\begin{multline}\label{eq:Diretta}
\int_{-T}^{T} \bigg| \sum_{n=n_0}^{\infty}C_{n}e^{i\omega_{n} t}+\overline{C_{n}}e^{-i\overline{\omega_{n}}t}
+R_{n}e^{r_{n} t}+D_{n}e^{i\zeta_{n} t}+\overline{D_{n}}e^{-i\overline{\zeta_{n}}t}\bigg|^2\ dt
\\
+\int_{-T}^{T} \bigg| \sum_{n=n_0}^{\infty}d_nD_{n}e^{i\zeta_{n} t}+\overline{d_nD_{n}}e^{-i\overline{\zeta_{n}}t}
+c_nC_{n}e^{i\omega_{n} t}+\overline{c_nC_{n}}e^{-i\overline{\omega_{n}}t}\bigg|^2\ dt
\\
\le
c(T)\sum_{ n= n_0}^{\infty}
\Big(|C_{n}|^2+ |d_nD_{n}|^2\Big)
\,.
\end{multline}
\end{theorem}
\begin{Proof}
Since the function $k^*(t)$ is positive, for $n_0\in\N$ to be chosen later we have
\begin{multline*}
\int_{-\infty}^{\infty} k^*(t)
\bigg|\sum_{n=n_0}^{\infty}C_{n}e^{i\omega_{n} t}+\overline{C_{n}}e^{-i\overline{\omega_{n}}t}
+R_{n}e^{r_{n} t}+D_{n}e^{i\zeta_{n} t}+\overline{D_{n}}e^{-i\overline{\zeta_{n}}t}\bigg|^2\ dt
\\
\le
4\int_{-\infty}^{\infty} k^*(t)\bigg|\sum_{n=n_0}^\infty C_{n}e^{i\omega_{n}t}
+\overline{C_n}e^{-i\overline{\omega_n}t}\bigg|^2\ dt
+4\int_{-\infty}^{\infty} k^*(t)\bigg|\sum_{n=n_0}^\infty R_{n}e^{r_{n}t}\bigg|^2\ dt
\\
+4\int_{-\infty}^{\infty} k^*(t)\bigg|\sum_{n=n_0}^\infty D_{n}e^{i\zeta_{n} t}+\overline{D_{n}}e^{-i\overline{\zeta_{n}}t}\bigg|^2\ dt\,.
\end{multline*}
We can apply Proposition \ref{pr:Fnd} to the first term and to the third one and Proposition \ref{pr:Rnd} to the second term. Therefore,  fixed $\varepsilon\in (0,1)$ and 
$T>\frac{\pi}{\ga\sqrt{1-\varepsilon}}$ there exists $n_0=n_0(\varepsilon)\in\N$ such that, thanks to inequalities \eqref{eq:Fndir}--\eqref{le:U1R} and in view also of \eqref{eq:cndn}, we get 
\begin{multline}\label{eq:firstd}
\int_{-\infty}^{\infty} k^*(t)
\bigg|\sum_{n=n_0}^{\infty}C_{n}e^{i\omega_{n} t}+\overline{C_{n}}e^{-i\overline{\omega_{n}}t}
+R_{n}e^{r_{n} t}+D_{n}e^{i\zeta_{n} t}+\overline{D_{n}}e^{-i\overline{\zeta_{n}}t}\bigg|^2\ dt
\\
\le
c(T)
\sum_{n= n_0}^\infty
\Big(|C_{n}|^2+ |d_nD_{n}|^2\Big)
\,.
\end{multline}
Moreover, in a similar way applying again Proposition \ref{pr:Fnd} and taking into account \eqref{eq:cndn} we have
\begin{multline*}
\int_{-\infty}^{\infty} k^*(t)
\bigg| \sum_{n=n_0}^{\infty}d_nD_{n}e^{i\zeta_{n} t}+\overline{d_nD_{n}}e^{-i\overline{\zeta_{n}}t}
+c_nC_{n}e^{i\omega_{n} t}+\overline{c_nC_{n}}e^{-i\overline{\omega_{n}}t}\bigg|^2\ dt
\\
\le
c(T)
\sum_{n= n_0}^\infty
\Big( |d_nD_{n}|^2+|C_{n}|^2\Big)
\,.
\end{multline*}
Combining \eqref{eq:firstd} with the above inequality and recalling the notation \eqref{eq:notation}  yields
\begin{equation*}
\int_{-\infty}^{\infty} k^*(t)
\big(|u_1^{n_0}(t)|^2+|u_2^{n_0}(t)|^2\big)\ dt
\le
c(T)
\sum_{n= n_0}^\infty
\Big(|C_{n}|^2+ |d_nD_{n}|^2\Big)
\,.
\end{equation*}
Now, we can consider the last inequality with the function $k^*$ 
replaced by the analogous one relative to $2T$ instead of $T$.
So, taking into account (\ref{eq:kcos}), 
we get 
\begin{equation*}
\int_{-2T}^{2T} \cos \frac{\pi t}{4T}
\big(|u_1^{n_0}(t)|^2+|u_2^{n_0}(t)|^2\big)\ dt
\le
c(2T)
\sum_{n= n_0}^\infty
\Big(|C_{n}|^2+ |d_nD_{n}|^2\Big),
\end{equation*}
whence, thanks to
$
\cos \frac{\pi t}{4T}\ge \frac1{\sqrt2}
$
for $|t|\le T$,
it follows
\begin{equation*}
\int_{-T}^{T} \big(|u_1^{n_0}(t)|^2+|u_2^{n_0}(t)|^2\big)\ dt
\le
\sqrt2c(2T)\sum_{n= n_0}^\infty
\Big(|C_{n}|^2+ |d_nD_{n}|^2\Big)\,.
\end{equation*}
This completes the proof.
\end{Proof}
Based on the approach performed in \cite{Ha}, the next result states that we can recover the finite number of missing terms  in the inverse and direct estimates.
We omit the proof, because it may be proved in much the same way as Proposition 5.8 and Proposition 5.20 of \cite{LoretiSforza3}.
We advise the reader to keep in mind  formulas \eqref{eq:vsum1} and \eqref{eq:notation}.

\begin{proposition}\label{pr:haraux-inv}
Let $\{\om_n\}_{n\in\N}$, $\{r_n\}_{n\in\N}$  and
$\{\zeta_{n}\}_{n\in\N}$ be sequences of pairwise  distinct numbers
such that $\om_n\not= \zeta_m$, $\om_n\not=\overline{\zeta_m}$, 
$r_n\not= i\om_m$, $r_n\not= i\zeta_m$, $r_n\not=-\eta$, $\zeta_{n}\not=0$, for any $n\,,m\in\N$, and
\begin{equation}\label{eq:ha1}
\lim_{n\to\infty}|\om_n|=\lim_{n\to\infty}|\zeta_{n}|=+\infty\,.
\end{equation}
Assume that there exists $n_0\in\N$ such that 
\begin{equation*}
\int_{0}^{T} \big(|u_1^{n_0}(t)|^2+|u_2^{n_0}(t)|^2\big)\ dt
\asymp
\sum_{ n= n_0}^{\infty}
\Big(|C_{n}|^2+ |d_nD_{n}|^2\Big)\,.
\end{equation*} 
Then, for any sequences $\{C_n\}$, $\{R_n\}$, $\{D_n\}$ and $\mathcal {E}\in\R$
we have
\begin{equation}\label{eq:haraux-inv22}
\int_{0}^{T} \big(|u_1(t)|^2+|u_2(t)|^2\big)\ dt
\asymp
\sum_{ n= 1}^{\infty}
\Big(|C_{n}|^2+ |d_nD_{n}|^2\Big) +|\mathcal {E}|^2 \,.
\end{equation}
\end{proposition}

\subsection{Inverse and direct inequalities}

We recall that
\begin{equation*}
\begin{split}
u_1(t) &=\sum_{n=1}^{\infty}\Big(C_{n}e^{i\omega_{n} t}+\overline{C_{n}}e^{-i\overline{\omega_{n}}t}
+R_{n}e^{r_{n} t}+D_{n}e^{i\zeta_{n} t}+\overline{D_{n}}e^{-i\overline{\zeta_{n}}t}\Big)
\,,
\\
u_2(t) &=\sum_{n=1}^{\infty}
\Big(d_nD_{n}e^{i\zeta_{n} t}+\overline{d_nD_{n}}e^{-i\overline{\zeta_{n}}t}
+c_nC_{n}e^{i\omega_{n} t}+\overline{c_nC_{n}}e^{-i\overline{\omega_{n}}t}\Big)
+\mathcal {E} e^{-\eta t}
\,,
\end{split}
\end{equation*}
where
\begin{equation}\label{eq:mathcalE}
|\mathcal {E}|^2\le 
M \sum_{ n= 1}^{\infty}
\Big(|C_{n}|^2+ |d_nD_{n}|^2\Big),
\qquad
(M>0)
\,.
\end{equation}

\begin{theorem}\label{th:inv.ingham1}
Let $\{\om_n\}_{n\in\N}$, $\{r_n\}_{n\in\N}$  and
$\{\zeta_{n}\}_{n\in\N}$ be sequences of pairwise  distinct numbers
such that $\om_n\not= \zeta_m$, $\om_n\not=\overline{\zeta_m}$, 
$r_n\not= i\om_m$, $r_n\not= i\zeta_m$, $r_n\not=-\eta$, $\zeta_{n}\not=0$, for any $n\,,m\in\N$.
Assume that  there exist
$\gamma>0$, $\alpha,\chi\in\R$, $n'\in\N$, $\mu>0$, $\nu> 1/2$, 
such that
\begin{equation*}
\liminf_{n\to\infty}({\Re}\om_{n+1}-{\Re}\om_{n})=\liminf_{n\to\infty}({\Re}\zeta_{n+1}-{\Re} \zeta_{n})=\gamma\,,
\end{equation*}
\begin{equation*}
\begin{split}
\lim_{n\to\infty}{\Im}\om_n&=\alpha>0
\,,
\\
\lim_{n\to\infty}r_n&=\chi<0\,,
\\
\lim_{n\to\infty}\Im \zeta_{n}&=0\,,
\end{split}
\end{equation*}
\begin{equation*}
|d_n|\asymp|\zeta_n|
\,,
\qquad
|c_n|\le\frac{M}{|\omega_n|}\,,
\end{equation*}
\begin{equation*}
|R_n|\le \frac{\mu}{n^{\nu}}\Big(|C_{n}|^2+ |d_nD_{n}|^2\Big)^{1/2}\,\quad\forall\ n\ge n'\,,
\qquad
|R_n|\le \mu\Big(|C_{n}|^2+ |d_nD_{n}|^2\Big)^{1/2}\,\quad\forall\ n\le n'\,.
\end{equation*}
Then, for $\gamma>4\alpha$ and
$T>\frac{2\pi}{\sqrt{\gamma^2-16\alpha^2}}$ we have 
 \begin{equation}\label{eq:inv.ingham}
 \int_{0}^{T} \big(|u_1(t)|^2+|u_2(t)|^2\big)\ dt
 \asymp
\sum_{ n= 1}^{\infty}
\Big(|C_{n}|^2+ |d_nD_{n}|^2\Big) \,.
\end{equation}
\end{theorem}
\begin{Proof}
Since $T>\frac{2\pi}{\sqrt{\gamma^2-16\alpha^2}}$, there exists $0<\varepsilon<1$ such that
$T>\frac{2\pi}{\sqrt{\gamma^2(1-\varepsilon)-16\alpha^2(1+\varepsilon)}}$. Therefore, thanks to Theorems \ref{th:gamma>4alpha} and \ref{th:Diretta}  we are able to employ Proposition \ref{pr:haraux-inv} obtaining
 \begin{equation*}
 \int_{0}^{T} \big(|u_1(t)|^2+|u_2(t)|^2\big)\ dt
 \asymp
\sum_{ n= 1}^{\infty}
\Big(|C_{n}|^2+ |d_nD_{n}|^2\Big) +|\mathcal {E}|^2 \,.
\end{equation*}
Finally, by \eqref{eq:mathcalE} we can get rid of the term $|\mathcal {E}|^2$ in the previous estimates,
and hence the proof is complete.
\end{Proof}

If we assume the condition $|C_n|\le M |d_nD_{n}|$ on the coefficients of the series
 instead of  $\gamma>4\alpha$, then we can make use of Theorem \ref{th:extracoe} instead of Theorem \ref{th:gamma>4alpha}, obtaining the observability inequalities with a better estimate for the control time: $T>\frac{2\pi}{\gamma}$. Precisely, the following result holds.

\begin{theorem}\label{th:inv.ingham11}
Let  assume the hypotheses of Theorem \ref{th:inv.ingham1} and the condition 
\begin{equation}
|C_n|\le M |d_nD_{n}|
\,.
\end{equation}
Then, for 
$T>\frac{2\pi}{\gamma}$ we have 
 \begin{equation}\label{eq:inv.ingham11}
 \int_{0}^{T} \big(|u_1(t)|^2+|u_2(t)|^2\big)\ dt
 \asymp
\sum_{ n= 1}^{\infty}
\Big(|C_{n}|^2+ |d_nD_{n}|^2\Big) \,.
\end{equation}

\end{theorem}

\section{Reachability results}

This section will be devoted to the proof of some
reachability results for  wave--wave coupled systems with a memory term.
\begin{theorem}\label{th:reachres}
Let $\beta<1/2$. For  any $T>\frac{2\pi}{\sqrt{1-4\beta^2}}$ and
$
(u_{i}^{0},u_{i}^{1})\in  L^{2}(0,\pi)\times H^{-1}(0,\pi)
$,
$i=1,2$,
 there exist $g_i\in L^2(0,T)$, $i=1,2$, such that the weak solution  $(u_1,u_2)$ of system 
\begin{equation}\label{eq:problem-usix}
\begin{cases}
\displaystyle 
u_{1tt}(t,x) -u_{1xx}(t,x)+\beta\int_0^t\ e^{-\eta(t-s)} u_{1xx}(s,x)ds+au_2(t,x)= 0\,,
\\
\phantom{u_{1tt}(t,x) -u_{1xx}(t,x)+\int_0^t\ k(t-s) u_{1xx}(s,x)ds+}
t\in (0,T)\,,\,\,\, x\in(0,\pi)
\\
\displaystyle
u_{2tt}(t,x) -u_{2xx}(t,x)+bu_1(t,x)= 0
\,,
\end{cases}
\end{equation} 
with boundary conditions
\begin{equation}\label{eq:bound-u1r}
u_1(t,0)=u_2(t,0)=0\,,\quad u_1(t,\pi)=g_1(t)\,,\quad u_2(t,\pi)=g_2(t)\qquad t\in (0,T) 
 \,,
\end{equation}
and null initial values 
\begin{equation}
u_i(0,x)=u_{it}(0,x)=0\qquad  x\in(0,\pi)\,,\quad i=1,2,
\end{equation} 
verifies the final conditions
\begin{equation}\label{eq:findataT}
u_i(T,x)=u_{i}^{0}(x)\,,\quad u_{it}(T,x)=u_{i}^{1}(x)\,,
\quad x\in(0,\pi),
\qquad i=1,2\,.
\end{equation}
\end{theorem}
\begin{Proof}
To prove our statement, we will apply the Hilbert Uniqueness Method described in Section \ref{se:HUM}.
Let 
$
H= L^2(0,\pi )
$
be endowed with the usual scalar product and norm 
$$
\|u\|_{L^2}:=\left(\int_0^\pi |u(x)|^{2}\ dx\right)^{1/2}\qquad
u\in L^2(0,\pi)\,.
$$
We consider the operator $L:D(L)\subset H\to H$ defined by $Lu=\displaystyle -u_{xx}$ for $u\in D(L):=H^2(0,\pi )\cap H_0^1(0,\pi )$.
It is well known that $L$ is a self-adjoint positive
 operator on $H$ with dense domain $D(L)$ and 
 $$D(\sqrt L)=H_0^1(0,\pi ).$$
 Moreover, $\{n^2\}_{n\ge1}$ is the sequence  of  eigenvalues for  $L$ and 
  $\{\sin(nx)\}_{n\ge1}$ is the sequence of the corresponding eigenvectors.
We can apply our spectral analysis, see Section \ref{se:specan}, to the adjoint system of (\ref{eq:problem-usix}) given by 
\begin{equation}\label{eq:adjointr}
\begin{cases}
\displaystyle 
z_{1tt}(t,x) -z_{1xx}(t,x)+\int_t^T\ k(s-t) z_{1xx}(s,x)ds+bz_2(t,x)= 0\,,\\
\hskip7cm
t\in (0,T)\,,\ x\in(0,\pi)
\\
\displaystyle
z_{2tt}(t,x) -z_{2xx}(t,x)+az_1(t,x)= 0\,,
\\
z_i(t,0)=z_i(t,\pi)=0\quad t\in [0,T]\,, 
\\
\hskip5cm  i=1,2,
\\
z_i(T,\cdot)=z_{i}^0\,,\quad z_{it}(T,\cdot)=z_{i}^{1}\,,
\end{cases}
\end{equation}
where the final data exhibit the following expansion in the basis $\{\sin(nx)\}_{n\ge1}$
\begin{equation*}
z_{i}^{0}(x)=\sum_{n=1}^{\infty}\alpha_{in}\sin(nx)\,,\quad
              z_{i}^{1}(x)=\sum_{n=1}^{\infty}\rho_{in}\sin(nx)
\,,\qquad  i=1,2\,.
\end{equation*} 
If we take 
$
(z_{i}^{0},z_{i}^{1})\in
H^1_0(0,\pi)\times L^2(0,\pi)
$,
$i=1,2$,
then one has
 \begin{equation}\label{eq:norms}
              \|z_{i}^{0}\|^2_{H^1_0}=\sum_{n=1}^\infty\alpha^2_{in} n^2,
              \quad
              \|z_{i}^{1}\|^2_{L^2}= \sum_{n=1}^\infty\rho^2_{in}
\,,\qquad i=1,2.
\end{equation} 
The backward system \eqref{eq:adjointr} is equivalent to the forward system
\begin{equation}\label{eq:forward}
\begin{cases}
\displaystyle 
u_{1tt}(t,x) -u_{1xx}(t,x)+\int_0^t\ k(t-s) u_{1xx}(s,x)ds+bu_2(t,x)= 0\,,\\
\hskip7cm
t\in (0,T)\,,\ x\in(0,\pi)
\\
\displaystyle
u_{2tt}(t,x) -u_{2xx}(t,x)+au_1(t,x)= 0\,,
\\
u_i(t,0)=u_i(t,\pi)=0\quad t\in [0,T]\,, 
\\
\hskip5cm  i=1,2,
\\
u_i(0,\cdot)=z_{i}^{0}\,,\quad u_{it}(0,\cdot)=z_{i}^{1}\,,
\end{cases}
\end{equation}
that is, if  $(u_1,u _2)$ is the solution of \eqref{eq:forward}, then the solution $(z_1,z_2)$ of \eqref{eq:adjointr} is given by 
\begin{equation*}
z_1(t,x)=u_1(T-t,x), \qquad z_2(t,x)=u_2(T-t,x)
\,.
\end{equation*}
Therefore, thanks to the representation for the solution of \eqref{eq:forward}, see Theorem \ref{th:repres}, we
can write $(z_1,z_2)$  in the following way, for any $(t,x)\in [0,T]\times [0,\pi]$
\begin{equation*}
z_1(t,x)=
\sum_{n=1}^{\infty}\Big(C_ne^{i\om_n(T-t)}
+\overline{C_n}e^{-i\overline{\om_n}(T-t)}+R_ne^{r_n(T-t)}+D_ne^{i\zeta_n(T-t)}
+\overline{D_n}e^{-i\overline{\zeta_n}(T-t)}\Big)\sin(n x)
\,,
\end{equation*}
\begin{multline*}
z_2(t,x)=\sum_{n=1}^{\infty}
\Big(d_nD_{n}e^{i\zeta_{n} (T-t)}+\overline{d_nD_{n}}e^{-i\overline{\zeta_{n}}(T-t)}
+c_nC_{n}e^{i\omega_{n} (T-t)}+\overline{c_nC_{n}}e^{-i\overline{\omega_{n}}(T-t)}\Big)\sin(n x)
\\
+e^{-\eta (T-t)}\sum_{n=1}^{\infty}E_n\sin(n x)
\,.
\end{multline*}
In particular, thanks also to \eqref{eq:norms}
we get
\begin{equation}\label{eq:equivn}
 \sum_{n= 1}^\infty n^2\Big( |C_n|^2+|d_nD_{n}|^2 \Big)
  \asymp
 \|z_{1}^{0}\|^2_{H^1_0}+\|z_{1}^{1}\|^2_{L^2}
 +\|z_{2}^{0}\|^2_{H^1_0}+\|z_{2}^{1}\|_{L^2}^2
 \,.
 \end{equation}
Moreover, for any $t\in [0,T]$
\begin{equation*}
z_{1x}(t,\pi)=
\sum_{n=1}^{\infty}(-1)^n n\Big(C_ne^{i\om_n(T-t)}
+\overline{C_n}e^{-i\overline{\om_n}(T-t)}+R_ne^{r_n(T-t)}+D_ne^{i\zeta_n(T-t)}
+\overline{D_n}e^{-i\overline{\zeta_n}(T-t)}\Big)
\,,
\end{equation*}
\begin{multline*}
z_{2x}(t,\pi)=\sum_{n=1}^{\infty}(-1)^n n
\Big(d_nD_{n}e^{i\zeta_{n} (T-t)}+\overline{d_nD_{n}}e^{-i\overline{\zeta_{n}}(T-t)}
+c_nC_{n}e^{i\omega_{n} (T-t)}+\overline{c_nC_{n}}e^{-i\overline{\omega_{n}}(T-t)}\Big)
\\
+e^{-\eta (T-t)}\sum_{n=1}^{\infty}(-1)^n n\, E_n
\,.
\end{multline*}
We can apply Theorem \ref{th:inv.ingham1}  to $(z_{1x}(t,\pi),z_{2x}(t,\pi))$.
Indeed, thanks to the above expressions for $z_{ix}(t,\pi)$, $i=1,2$, and  \eqref{eq:inv.ingham} we have   
\begin{equation*}
 \int_{0}^{T}  \big(|z_{1x}(t,\pi)|^2+|z_{2x}(t,\pi)|^2\big)\ dt
  \asymp
 \sum_{n= 1}^\infty n^2\Big( |C_n|^2+|d_nD_{n}|^2 \Big)\,,
 \end{equation*}
and hence by \eqref{eq:equivn}
we get
\begin{equation}\label{eq:obses}
 \int_{0}^{T}  \big(|z_{1x}(t,\pi)|^2+|z_{2x}(t,\pi)|^2\big)\ dt
  \asymp
 \|z_{1}^{0}\|^2_{H^1_0}+\|z_{1}^{1}\|^2_{L^2}
 +\|z_{2}^{0}\|^2_{H^1_0}+\|z_{2}^{1}\|_{L^2}^2
 \,.
 \end{equation}
 Therefore, we have proved Theorem \ref{th:uniqueness}. Furthermore, 
we consider the linear operator $\Psi$   introduced in Section \ref{se:HUM} and, thanks to \eqref{eq:psi0},  defined by 
\begin{equation*}
\Psi(z_{1}^{0},z_{1}^{1},z_{2}^{0},z_{2}^{1})=(-u_{1t}(T,\cdot),u_{1}(T,\cdot),-u_{2t}(T,\cdot),u_{2}(T,\cdot))
\,,
\end{equation*}
where $(u_{1},u_{2})$ is the weak solution of system \eqref{eq:problem-usix}. We have that
$$\Psi:H^1_0(0,\pi)\times L^2(0,\pi)\times H^1_0(0,\pi)\times L^2(0,\pi)\to H^{-1}(0,\pi)\times L^2(0,\pi)\times H^{-1}(0,\pi)\times L^2(0,\pi)$$
 is an isomorphism.
Therefore, for $
(u_{i}^{0},u_{i}^{1})\in  L^{2}(0,\pi)\times H^{-1}(0,\pi)
$,
$i=1,2$, there exists one and only one $(z_{1}^{0},z_{1}^{1},z_{2}^{0},z_{2}^{1})\in H^1_0(0,\pi)\times L^2(0,\pi)\times H^1_0(0,\pi)\times L^2(0,\pi)$
such that 
\begin{equation*}
\Psi(z_{1}^{0},z_{1}^{1},z_{2}^{0},z_{2}^{1})
=(-u_{1}^{1},u_{1}^{0},-u_{2}^{1},u_{2}^{0})
\,.
\end{equation*}
Finally, if we consider the solution $(z_{1},z_{2})$ of system \eqref{eq:adjointr} with final data given by the unique $(z_{1}^{0},z_{1}^{1},z_{2}^{0},z_{2}^{1})$, then
the control functions required by the statement are given by
\begin{equation*}
g_1(t)=z_{1x}(t,\pi)-\beta\int_t^T\ e^{-\eta(s-t)}z_{1x}(s,\pi)ds\,,
\qquad
g_2(t)=z_{2x}(t,\pi)\,,
\end{equation*}
that is, our proof is complete.
\end{Proof}

\end{document}